\RequirePackage{fix-cm}
\documentclass[smallextended]{svjour3}
\smartqed

\usepackage{amsmath,amssymb,mathtools}
\usepackage{bm}
\usepackage{booktabs}
\usepackage{microtype}
\usepackage{xcolor}
\usepackage{tikz}
\usetikzlibrary{arrows.meta}
\usepackage{graphicx}
\usepackage{comment}
\usepackage{algorithm}
\usepackage{algorithmic}
\usepackage{enumitem}
\usepackage{multirow}
\usepackage{longtable}
\usepackage[hidelinks]{hyperref}
\usepackage{hyperref}

\newcommand{\Ssym}{\mathbb{S}}
\newcommand{\one}{\mathbf{1}}
\newcommand{\zero}{\mathbf{0}}
\newcommand{\I}{\mathbf{I}}
\newcommand{\J}{\mathbf{J}}
\newcommand{\Diag}{\operatorname{Diag}}
\newcommand{\diag}{\operatorname{diag}}
\newcommand{\tr}{\operatorname{tr}}

\newcommand{\ip}[2]{\left\langle #1,#2\right\rangle}
\newcommand{\inprod}[2]{\left\langle #1,#2\right\rangle}
\newcommand{\norm}[1]{\left\lVert #1\right\rVert}

\newcommand{\noise}{\sigma}

\begin{document}


\title{Optimality Guarantees for Minimum Variance Subset Selection via Semidefinite Programming}
\titlerunning{SDP tightness and exact solution of MV$k$SP}

\author{Antonio M. Sudoso \and Daniel Aloise}
\authorrunning{A. M. Sudoso and D. Aloise}

\institute{A. M. Sudoso \at
Department of Mathematics ``Tullio Levi-Civita'', University of Padua,
Via Trieste 63, 35131 Padova, Italy \\
\email{antoniomaria.sudoso@unipd.it} \\
ORCID: 0000-0002-2936-9931
\and
D. Aloise \at
Department of Computer and Software Engineering, Polytechnique Montr\'eal,
2500 chemin de Polytechnique, Montr\'eal, QC H3T 1J4, Canada \\
\email{daniel.aloise@polymtl.ca} \\
ORCID: 0000-0002-9876-2921}

\date{}
\maketitle

\begin{abstract}
Given a set of data points in Euclidean space and a prescribed
cardinality $k$, the minimum variance $k$-subset problem asks for $k$ points
minimizing their average squared distance from their centroid. We develop a
semidefinite programming (SDP) framework for exact recovery and global
optimization. First, we show that the natural SDP
relaxation fails Slater's condition. Identifying the proper face containing
its feasible region yields an equivalent, lower-dimensional, strictly
feasible formulation that enjoys strong duality. A primal--dual characterization of
optimality then gives geometric necessary and sufficient conditions
under which the SDP relaxation is guaranteed to be tight. These conditions reveal how tightness depends on
the internal geometry of the selected subset and on its separation from the
remaining points.
Second, we embed the SDP relaxation in a branch-and-cut
algorithm combining first-order optimization, valid inequalities and dimension-reducing branching. Extensive experiments on public data
sets, with up to $1{,}500$ points and different values of $k$, show that the
method strongly outperforms general-purpose exact solvers as the number of points and cardinality increase.
\end{abstract}




\section{Introduction}
\label{sec:introduction}

Given a set of $n$ points in $\mathbb{R}^d$ and an integer $1<k<n$, the minimum variance $k$-subset problem (MV$k$SP) asks for $k$ points minimizing the sum of their squared Euclidean distances to their centroid. An optimal subset is thus the group of $k$ points of minimum variance, hence the most cohesive one, and its centroid a representative estimated from this group. Both aspects have been studied. Finding $k$ points that form a small set under a given measure, the variance among others, is a classical problem in computational geometry~\cite{aggarwal1991finding}. In robust statistics, the centroid of an optimal subset is the impartial trimmed mean~\cite{gordaliza1991breakdown}, a location estimate insensitive to outliers obtained by discarding the $n-k$ observations whose removal minimizes the variance of the remaining ones. These two aspects meet in the analysis of measurement data, where a known number of identical objects is observed under Gaussian noise among unrelated ones. There, an optimal subset identifies these objects, and its centroid is the maximum-likelihood estimate of their common characteristics~\cite{KelmanovPyatkin2011}.

From an optimization perspective, the MV$k$SP can be formulated as a binary quadratic program with a cardinality constraint, and is closely related to several well-studied problems in the literature. In particular, it admits an equivalent formulation as the problem of finding the heaviest $k$-subgraph~\cite{billionnet2017using,pereira2026branch,malick2012solving,sotirov2020solving} in a complete graph whose edge weights are given by the negative squared Euclidean distances between the data points. The MV$k$SP is also related to the maximum diversity problem~\cite{duarte2007tabu,marti2010branch,martinez2026stochastic}, which seeks to select $k$ points that maximize the sum of pairwise distances among them, in contrast to the present problem, which minimizes variance. It is further connected to the minimum sum-of-squares clustering (MSSC) problem~\cite{aloise2012improved,piccialli2025global,sudoso2025column}, which partitions all data points into clusters and minimizes the sum of within-cluster squared deviations. Unlike MSSC, however, the MV$k$SP focuses on selecting a single subset of fixed cardinality, without requiring a partition of the entire dataset.

The decision version of the MV$k$SP is strongly NP-complete when both $k$ and $d$ are part of the input~\cite{KelmanovPyatkin2011}. Its mathematical programming formulation as a binary quadratic program with a cardinality constraint leads to a nonconvex optimization problem whose standard linear relaxations provide weak bounds and lead to large branch-and-bound trees~\cite{malick2012solving,pereira2026branch}. These challenges motivate the development of stronger formulations and bounding techniques that better exploit the quadratic structure of the problem. In particular, semidefinite programming (SDP) relaxations have proven effective for a variety of quadratic combinatorial optimization problems, including max-cut~\cite{hrga2021madam,gusmeroli2022biqbin} and clique-related problems~\cite{ames2014guaranteed,sudoso2025semidefinite}, by providing strong bounds that can be embedded in branch-and-bound schemes to tackle large instances.

\paragraph{Contributions.}
We study the MV$k$SP from both theoretical and algorithmic perspectives, with the goal of solving the problem to global optimality. Our main contributions are as follows.

\begin{itemize}
    \item We show that the natural SDP relaxation fails Slater's condition because its affine constraints restrict feasibility to a proper face of the positive semidefinite cone. We explicitly identify this face and perform one facial reduction step, obtaining an equivalent, smaller SDP relaxation that is strictly feasible and satisfies strong duality.

    \item Using strong duality, we derive optimality conditions and translate them into geometric sufficient and necessary conditions for tightness. For a candidate subset whose points are equidistant from its centroid, we prove that the SDP relaxation recovers it when its diameter is smaller than both its separation from the unselected points and the separation among the unselected points. For arbitrary subsets, we introduce an internal shape parameter and an adjusted separation measure that together yield a polynomial-time checkable sufficient condition for exact recovery. We also derive a necessary condition requiring every unselected point to lie sufficiently far from the centroid.

    \item We validate the theory on synthetic instances with progressively perturbed planted solutions. The sufficient condition certifies exact recovery under substantial perturbations, while the necessary condition provides a useful indication of lost planted optimality.

    \item We exploit the facially reduced SDP relaxation algorithmically by embedding it within a branch-and-cut framework. To handle large-scale instances, the SDP relaxation is approximately solved by a first-order method, together with a post-processing procedure that provides valid lower bounds. These bounds are further strengthened by separating valid inequalities and adding them through a cutting-plane scheme. Branching produces child relaxations of the same form and smaller dimension.

    \item We benchmark the exact method on datasets with up to $1{,}500$ points across a wide range of dimensions and cardinality parameters. The proposed approach is compared with state-of-the-art general-purpose global optimization solvers and substantially outperforms them on the most challenging instances. 
\end{itemize}

The remainder of the paper is organized as follows. Section~\ref{sec:formulation} gives
the MV$k$SP formulation and positions it within the $k$-subgraph literature.
Section~\ref{sec:sdp} derives and facially reduces the SDP relaxation. Section~\ref{sec:tightness} develops
the exact-recovery conditions. Section~\ref{sec:algorithm} presents the branch-and-cut
algorithm, and Section~\ref{sec:experiments} reports the computational study. Section~\ref{sec:concl} concludes
the paper.

\section{Problem Formulation}
\label{sec:formulation}

Let $\mathcal{P} = \{p_1,\dots,p_n\} \subset \mathbb{R}^d$ be a finite set of  points. 
Given an integer $1 < k < n$, our goal is to select a subset of exactly $k$ points whose empirical variance about its centroid is minimized.
We introduce binary decision variables $x_i \in \{0,1\}$ for $i \in \{1, \dots, n\}$, where $x_i = 1$ indicates that point $p_i$ is selected and $x_i = 0$ otherwise. 
The centroid of the selected points is thus given by
\begin{equation*}
\mu(x) = \frac{1}{k}\sum_{i=1}^n x_i p_i.
\end{equation*}
The corresponding empirical variance of the selected subset is
\begin{equation*}
\frac{1}{k}\sum_{i=1}^n x_i \|p_i - \mu(x)\|^2,
\end{equation*}
which leads to the combinatorial optimization problem
\begin{equation}
\label{eq:original}
\begin{aligned}
\min_{x} \quad
& \frac{1}{k}\sum_{i=1}^n x_i \|p_i - \mu(x)\|^2 \\
\text{s.t.}\quad 
& \sum_{i=1}^n x_i = k, \\
& x_i\in\{0,1\}, \qquad \forall i \in \{1,\ldots,n\}.
\end{aligned}
\end{equation}
By Huygens' theorem (also known as the centroid identity  \cite{hastie2009elements}), the within-cluster sum of squared distances in the objective function can be expressed in terms of pairwise squared distances as
\begin{align*}
    \frac{1}{k}\sum_{i=1}^n x_i \lVert p_i - \mu(x)\rVert^2
    =
    \frac{1}{2k^2}
    \sum_{i=1}^n \sum_{j=1}^n
    x_i x_j \lVert p_i-p_j\rVert^2.
\end{align*}
Hence, Problem \eqref{eq:original} can be equivalently formulated as
\begin{equation}
\label{eq:quadratic_distance}
\begin{aligned}
\min_{x} \quad
& \frac{1}{2k^2} \sum_{i=1}^n \sum_{j=1}^n x_i x_j \|p_i - p_j \|^2 \\
\text{s.t.}\quad
& \sum_{i=1}^n x_i = k, \\
& x_i\in\{0,1\}, \qquad \forall i \in \{1,\ldots,n\}.
\end{aligned}
\end{equation}
Thus, Problem \eqref{eq:quadratic_distance} is a nonconvex quadratic binary optimization problem with a linear constraint.

The MV$k$SP belongs to the class of cardinality-constrained subset selection problems. It is a direct instance of the {heaviest $k$-subgraph problem} (H$k$SP)~\cite{billionnet2005different}, which asks for a $k$-vertex subset of an edge-weighted graph maximizing the sum of its induced edge weights. The unweighted special case of H$k$SP, in which all edge weights are restricted to $\{0,1\}$, is the more extensively studied {densest $k$-subgraph problem} (D$k$SP)~\cite{sotirov2020solving}.
To establish this connection, consider the complete graph on the points in $\mathcal{P}$ and define $w_{ij}=-\|p_i-p_j\|^2$ for $i\neq j$, with $w_{ii}=0$. Since $k$ is fixed, Problem~\eqref{eq:quadratic_distance} is equivalent to maximizing $\frac{1}{2}\sum_{i=1}^n\sum_{j=1}^n w_{ij}x_i x_j$ under the same cardinality and binary constraints, which is precisely H$k$SP.
The nonpositivity of the weights can be removed without changing the set of optimal solutions. Specifically, let $M\geq\max_{i<j}\|p_i-p_j\|^2$ and define $\widetilde w_{ij}=M-\|p_i-p_j\|^2$ for $i\neq j$, with $\widetilde w_{ii}=0$. For every feasible $x$, this transformation increases the maximization objective by the constant $M\binom{k}{2}$, yielding nonnegative edge weights while preserving the set of optimal solutions. 

Several exact approaches have been developed for H$k$SP~\cite{sotirov2020solving}. Billionnet~\cite{billionnet2005different} proposed mixed-integer linear programming formulations, while Billionnet et al.~\cite{billionnet2009improving} used semidefinite programming to construct an exact convex quadratic reformulation subsequently solved by branch-and-bound. Malick and Roupin~\cite{malick2012solving} and Krislock et al.~\cite{krislock2016computational} developed semidefinite branch-and-bound methods that solved instances with approximately 120--160 vertices. More recently, Pereira and da Cunha~\cite{pereira2026branch} proposed SDP-based branch-and-bound algorithms with bounds computed through Lagrangian relaxation. These studies indicate that exact methods are effective on small and medium-sized instances but remain computationally challenging at larger scales, particularly for dense graphs.
Moreover, some exact H$k$SP solvers require integer coefficients and therefore cannot be applied directly to general MV$k$SP instances.

Beyond this connection to existing exact methods, MV$k$SP has structure that motivates a specialized approach. Generic H$k$SP and D$k$SP instances may be defined on sparse or dense graphs, whereas MV$k$SP induces a complete graph with weights derived from squared Euclidean distances. This precludes exploiting graph sparsity but provides geometric structure that specialized algorithms can exploit. Moreover, while D$k$SP seeks a subset inducing many edges, MV$k$SP seeks one with small pairwise squared distances, reflecting geometric compactness rather than combinatorial density.

\subsection*{Notation}
Throughout this work, $\mathbb{S}^n$ denotes the set of $n \times n$ real symmetric matrices, {$\mathbb{S}_+^n$} the cone of positive semidefinite matrices, and $\mathbb{R}^n$ and $\mathbb{R}^{m \times n}$ the spaces of real vectors and matrices, respectively. For $X \in \mathbb{S}^n$, $X \succeq 0$ means that $X$ is positive semidefinite. Vectors $\zero_n$ and $\one_n$ are the all-zeros and all-ones vectors in $\mathbb{R}^n$, with subscripts omitted when the dimension is clear, while $\mathbf{I}_n$ and $\mathbf{J}_n$ are the identity and all-ones matrices in $\mathbb{S}^n$. For $x \in \mathbb{R}^n$, $x^\perp := \{y \in \mathbb{R}^n : x^\top y = 0\}$ denotes its orthogonal complement. Given an index set $S$, $X_{SS} := (X_{ij})_{i,j\in S}$ is the principal submatrix of $X$, and $x_S := (x_i)_{i\in S}$ the corresponding subvector of $x$. Moreover, $\mathrm{Diag}(x)$ is the diagonal matrix with diagonal $x$, whereas $\mathrm{diag}(X)$ is the vector of diagonal entries of $X$. The Frobenius inner product is $\langle A,B\rangle := \mathrm{tr}(B^\top A)$, with associated norm $\|A\|_F := \sqrt{\langle A,A\rangle}$. Finally, $A \circ B$ denotes the Hadamard product, defined by $(A\circ B)_{ij} := A_{ij}B_{ij}$.

\section{Semidefinite programming relaxation}
\label{sec:sdp}
Let $D$ be the squared Euclidean distance matrix with entries $D_{ij} = \|p_i-p_j\|^2$ for all $i, j \in \{1, \dots, n\}$.
To derive a tractable convex relaxation, we lift the quadratic term by introducing
a matrix variable $X \in \mathbb{S}^n$ intended to satisfy $X = xx^\top$.
Using $x^\top D x = \langle D, X\rangle$, problem
\eqref{eq:quadratic_distance} admits the lifted representation
\begin{equation}
\begin{aligned}
\min_{x,X} \quad
& \frac{1}{2k^2} \langle D, X\rangle \\
\text{s.t.}\quad
& \mathbf{1}_n^\top x = k,\\
& X = xx^\top,\\
& x \in \{0,1\}^n.
\end{aligned}
\end{equation}
Following~\cite{malick2012solving}, we strengthen the lifted formulation by exploiting the cardinality constraint. Specifically, under $X=xx^\top$, the equality $\mathbf{1}_n^\top x=k$ implies
$X\mathbf{1}_n=xx^\top\mathbf{1}_n=x(\mathbf{1}_n^\top x)=kx$, while binary feasibility yields $\operatorname{diag}(X)=x$, $x \geq \zero_n$, and $X \geq \zero_{n \times n}$ elementwise.
Furthermore, defining the lifted matrix
\[
Y = \begin{pmatrix}
X & x \\
x^\top & 1
\end{pmatrix},
\]
the constraint $X = xx^\top$ is equivalent to $Y \succeq 0$ and $\operatorname{rank}(Y)=1$.
Therefore, retaining the above valid linear and nonnegativity constraints and dropping the rank-one constraint yields the doubly nonnegative relaxation
\begin{equation}
\label{eq:sdp_relaxation}
\begin{aligned}
\min_{Y \in \mathbb{S}^{n+1}} \quad
& \frac{1}{2k^2} \langle D, X \rangle \\
\text{s.t.}\quad
& \mathbf{1}_n^\top x = k, \\
& \operatorname{diag}(X) = x, \\
& X\mathbf{1}_n = kx, \\ 
& Y =
\begin{pmatrix}
X & x \\
x^\top & 1
\end{pmatrix} \succeq 0, \ Y \geq \mathbf{0}_{(n+1) \times (n+1)}.
\end{aligned}
\end{equation}
Problem~\eqref{eq:sdp_relaxation} is convex and can be solved by semidefinite programming tools.
Moreover, every feasible binary solution $x$ of problem \eqref{eq:quadratic_distance} induces a feasible
point $(x, X)$ of \eqref{eq:sdp_relaxation} by setting $X = xx^\top$. Hence,
problem \eqref{eq:sdp_relaxation} provides a valid lower bound on the optimal value of the
original problem.

The relaxation \eqref{eq:sdp_relaxation} does not satisfy Slater's condition, i.e., it admits no strictly feasible
solution. We now perform one step of facial reduction to obtain an equivalent SDP
relaxation that is strictly feasible \cite{borwein1981regularizing}. For an introduction to facial reduction, we refer the interested reader to~\cite{drusvyatskiy2017many} and the references therein.
Indeed, the affine constraints force all feasible matrices $Y$ to lie in a proper face of the cone $\mathbb{S}_+^{n+1}$.
To see this, consider the vector
\[
v =
\begin{pmatrix}
\mathbf{1}_n \\
- k
\end{pmatrix} \in \mathbb{R}^{n+1}.
\]
For every feasible $Y$, the constraints $X\mathbf{1}_n = kx$ and $\mathbf{1}_n^\top x = k$ imply
\[
Y v =
\begin{pmatrix}
X\mathbf{1}_n - kx \\
x^\top \mathbf{1}_n - k
\end{pmatrix}
=
\begin{pmatrix}
\mathbf{0}_n \\
0
\end{pmatrix}.
\]
Hence, every feasible matrix $Y$ is singular and therefore no feasible point satisfies $Y \succ 0$. Equivalently, with the exposing matrix $vv^\top \succeq 0$, we have {$\inprod{vv^\top}{Y}$} $= v^\top Y v = 0$ for every feasible $Y$. Thus, the feasible set is contained in the proper exposed face
\[
\mathcal{F} = \left\{ Y \in \mathbb{S}_+^{n+1} : Yv = \mathbf{0}_{n+1} \right\} = \left\{  Y \in \mathbb{S}_+^{n+1} : \inprod{vv^\top}{Y} = 0\right\}.
\]
To parametrize this face, let $V \in \mathbb{R}^{(n+1)\times n}$ be a matrix whose columns form a basis of the orthogonal complement of $v$, i.e., $\textrm{rank}(V) = n$ and $V^\top v=\mathbf{0}_n$. A convenient choice is
\[
V =
\begin{pmatrix}
\mathbf{I}_n \\
\frac{1}{k}\mathbf{1}_n^\top
\end{pmatrix}.
\]
Since $\textrm{range}(V) = v^\perp$, every $Y\in\mathcal{F}$ admits a representation $Y = V Z V^\top$ for some $Z \in \mathbb{S}_+^n$. We thus obtain an equivalent reformulation of \eqref{eq:sdp_relaxation} in the variable $Z$ expressed as the following facially reduced semidefinite program:
\begin{equation}
\label{eq:fr_sdp}
\begin{aligned}
\min_{Z \in \mathbb{S}^n} \quad 
& \frac{1}{2k^2} \langle D, Z \rangle \\
\text{s.t.}\quad 
& Z \mathbf{1}_n = k\operatorname{diag}(Z), \\
& \textrm{tr}(Z) = k, \\
& Z \succeq 0, \ Z \geq \mathbf{0}_{n \times n}.
\end{aligned}
\end{equation}
The equivalence is formally stated in the next proposition.

\begin{proposition}
\label{prop:fr-equivalence}
The SDP relaxations \eqref{eq:sdp_relaxation} and \eqref{eq:fr_sdp} are equivalent.
\end{proposition}
\begin{proof}
Let $Y$ be feasible for problem \eqref{eq:sdp_relaxation}.
As shown above, $Yv = \mathbf{0}_{n+1}$, hence $Y \in \mathcal{F}$. Therefore, there exists $Z \in \mathbb{S}_+^n$ such that $Y = V Z V^\top$. Expanding this factorization gives
\[
Y =
\begin{pmatrix}
Z & \frac{1}{k} Z \mathbf{1}_n \\
\frac{1}{k} \mathbf{1}_n^\top Z & \frac{1}{k^2} \mathbf{1}_n^\top Z \mathbf{1}_n
\end{pmatrix}.
\]
Matching blocks yields
\[
X = Z, \quad
x = \frac{1}{k} Z \mathbf{1}_n, \quad
\mathbf{1}_n^\top Z \mathbf{1}_n = k^2.
\]
Substituting into the constraints we obtain $Z \mathbf{1}_n = k \operatorname{diag}(Z)$ and $\mathbf{1}_n^\top Z \mathbf{1}_n = k^2$.
Note that under the constraint $Z\mathbf{1}_n = k \operatorname{diag}(Z)$, the identity $\mathbf{1}_n^\top Z \mathbf{1}_n = k \operatorname{tr}(Z)$ holds. Consequently, constraints $\operatorname{tr}(Z)=k$ and
$\mathbf{1}_n^\top Z \mathbf{1}_n=k^2$ are equivalent. 
Finally $Y \ge \mathbf{0}_{(n+1) \times (n+1)}$ implies $Z = X \ge \mathbf{0}_{n \times n}$, hence $Z$ is feasible for problem \eqref{eq:fr_sdp}.
The objective values coincide since
\begin{align*}
    \frac{1}{2k^2}\langle D,X\rangle = \frac{1}{2k^2} \langle D,Z \rangle.
\end{align*}
Conversely, let $Z$ be feasible for problem \eqref{eq:fr_sdp}, and set
\[
x = \frac{1}{k} Z \mathbf{1}_n, \quad X = Z, \quad
Y = \begin{pmatrix} X & x \\ x^\top & 1 \end{pmatrix}.
\]
Since $Z \mathbf{1}_n = k \operatorname{diag}(Z)$ and $\operatorname{tr}(Z) = k$, we have $\mathbf{1}_n^\top Z \mathbf{1}_n=k \operatorname{tr}(Z) = k^2$. Therefore
\[
Y =
\begin{pmatrix}
Z & \frac{1}{k} Z \mathbf{1}_n \\
\frac{1}{k} \mathbf{1}_n^\top Z & \frac{1}{k^2} \mathbf{1}_n^\top Z \mathbf{1}_n
\end{pmatrix}
= V Z V^\top,
\]
and thus $Y \succeq 0$. Also, since $Z \ge \mathbf{0}_{n \times n}$, we have $x \geq \mathbf{0}_n$ so it follows that $Y \ge \mathbf{0}_{(n+1) \times (n+1)}$. Moreover,
\[
\operatorname{diag}(X) = \operatorname{diag}(Z) = \frac{1}{k} Z \mathbf{1}_n = x,
\]
\[
X \mathbf{1}_n = Z \mathbf{1}_n = kx,
\]
\[
\mathbf{1}_n^\top x = \frac{1}{k} \mathbf{1}_n^\top Z \mathbf{1}_n = k.
\]
Hence $Y$ is feasible for problem \eqref{eq:sdp_relaxation}.
The objective values coincide because $X = Z$.
\qed
\end{proof}
 
The facially reduced formulation is advantageous from a numerical perspective, as it restores strict feasibility on a lower-dimensional cone, which can significantly improve numerical stability and computational efficiency for both first-order and interior-point methods \cite{hu2023facial}.
We henceforth work with relaxation~\eqref{eq:fr_sdp}. We verify
Slater's condition for both the primal and the dual, thereby establishing attainment and
strong duality. This provides the foundation for deriving
necessary and sufficient tightness conditions and translating them into
geometric conditions.

\noindent

Let \(u\in\mathbb{R}^n\) and \(\tau\in\mathbb{R}\) be the multipliers associated with the affine constraints \(Z\mathbf{1}_n=k\operatorname{diag}(Z)\) and \(\langle \mathbf{I}_n,Z\rangle=k\), respectively, and let \(\Gamma,Q\in\mathbb{S}^n\) be the multipliers associated with \(Z\ge\mathbf{0}_{n \times n}\) and \(Z\succeq0\), respectively.
 Using the standard derivation in Lagrangian duality theory \cite{boyd2004convex}, the dual problem of \eqref{eq:fr_sdp} is given by
\begin{equation}
\label{eq:fr_sdp_dual}
\begin{aligned}
\max_{u,\tau,\Gamma,Q}\quad
& k \tau \\
\text{s.t.}\quad
&
\Diag(u) - \frac{1}{2k}\left(u\mathbf{1}_n^\top + \mathbf{1}_n u^\top\right)
+
\tau \mathbf{I}_n
+
\Gamma
+
Q
=
\frac{1}{2k^2} D, \\
& Q\succeq0,\ \Gamma \ge\mathbf{0}_{n \times n}.
\end{aligned}
\end{equation}

\begin{proposition}
\label{prop:slater}
The primal problem \eqref{eq:fr_sdp} and the dual problem \eqref{eq:fr_sdp_dual} are strictly feasible. 
\end{proposition}

\begin{proof}
For the primal problem, let
\[
\bar Z=\alpha \textbf{I}_n+\beta \mathbf{J}_n,
\qquad
\alpha=\frac{k(n-k)}{n(n-1)},
\qquad
\beta=\frac{k(k-1)}{n(n-1)}.
\]
Since \(1<k<n\), one has
\(\alpha>0\) and \(\beta>0\). The eigenvalues of \(\bar Z\)
are \(\alpha\), with multiplicity \(n-1\), and
\(\alpha+n\beta=k^2/n\), with multiplicity one. Thus
\(\bar Z\succ0\). Moreover, since its diagonal entries are equal to \(\alpha+\beta>0\) and its off-diagonal entries are equal to \(\beta>0\), we also have \(\bar Z>\textbf{0}_{n \times n}\) entrywise.
Furthermore,
\[
\operatorname{diag}(\bar{Z})=(\alpha+\beta)\mathbf{1}_n
=
\left(\frac{k(n-k)+k(k-1)}{n(n-1)} \right)\mathbf{1}_n
=
\frac{k}{n}\mathbf{1}_n,
\]
and
\[
\bar{Z} \mathbf{1}_n
=
(\alpha+n\beta) \mathbf{1}_n
=
\left(
\frac{k(n-k)}{n(n-1)}
+
\frac{nk(k-1)}{n(n-1)}
\right)\mathbf{1}_n
=
\frac{k^2}{n}\mathbf{1}_n
=
k\,\operatorname{diag}(\bar{Z}).
\]
Finally,
\[
\operatorname{tr}(\bar{Z})
=
n\alpha+n\beta
=
\frac{k(n-k)+k(k-1)}{n-1}
=
k.
\]
Thus \(\bar{Z}\) satisfies the affine constraints of \eqref{eq:fr_sdp} and lies in the interior of both cones. Hence, \(\bar{Z}\) is strictly feasible for \eqref{eq:fr_sdp}.

For the dual problem, let $C:=D/(2k^2)$, choose $\varepsilon>0$,
and set $\bar u=\mathbf{0}_n$ and
$\bar\Gamma=\varepsilon\mathbf{J}_n$.
Choose $\bar\tau<\lambda_{\min}(C-\bar\Gamma)$ and define
$\bar Q:=C-\bar\tau \mathbf{I}_n-\bar\Gamma$. Then
\[
\lambda_{\min}(\bar Q)
=\lambda_{\min}(C-\bar\Gamma)-\bar\tau>0,
\]
so $\bar Q\succ0$. Since $\bar u=\mathbf{0}_n$ and
$\bar\tau I_n+\bar\Gamma+\bar Q=C$, the dual equality holds. Thus, the point \((\bar u,\bar\tau,\bar\Gamma,\bar Q)\)
satisfies the affine constraint of \eqref{eq:fr_sdp_dual} and lies in the interior of both dual cones. Therefore,
\((\bar u,\bar\tau,\bar\Gamma,\bar Q)\)
is strictly feasible for \eqref{eq:fr_sdp_dual}.
\end{proof}


Since $D$ and every feasible $Z$ are entrywise nonnegative,
the primal objective is bounded below by zero.
Together with the strict feasibility of both problems established
in Proposition~\ref{prop:slater}, Slater's theorem implies that
both optimal values are attained and there is no duality gap.
Consequently, the primal--dual optimality conditions are necessary
and sufficient for optimality~\cite{vandenberghe1996semidefinite}.

\begin{lemma}[Complementary slackness]
\label{lem:complementarity}
Let \(Z\) be primal feasible, and let $(u,\tau,\Gamma,Q)$ form a
dual-feasible point. They are primal--dual optimal if and only if
\(\inprod{\Gamma}{Z}=0\) and \(\inprod{Q}{Z}=0\).
\end{lemma}

The optimality conditions in Lemma~\ref{lem:complementarity} can equivalently be stated as \(\Gamma\circ Z=\mathbf{0}_{n \times n}\) and \(QZ=ZQ=\mathbf{0}_{n \times n}\). Indeed, since \(\Gamma,Z\geq\mathbf{0}_{n \times n}\) entrywise, \(\inprod{\Gamma}{Z}=0\) is equivalent to \(\Gamma\circ Z=\mathbf{0}_{n \times n}\). Moreover, since \(Q,Z\succeq 0\), we have \(0=\inprod{Q}{Z}=\operatorname{tr}(QZ)=\|Q^{1/2}Z^{1/2}\|_F^2\), which implies \(Q^{1/2}Z^{1/2}=\mathbf{0}_{n \times n}\), and hence \(QZ=\mathbf{0}_{n \times n}\). Since \(Q\) and \(Z\) are symmetric, \(ZQ=(QZ)^\top=\mathbf{0}_{n \times n}\). Therefore, \(\inprod{Q}{Z}=0\) is equivalent to \(QZ=ZQ=\mathbf{0}_{n \times n}\).

\section{Tightness of the SDP relaxation}
\label{sec:tightness}

We now investigate conditions under which the SDP
relaxation~\eqref{eq:fr_sdp} is tight, namely, when it admits a
rank-one optimal solution corresponding to an optimal solution of the
original combinatorial problem. More precisely, let
\(S\subseteq\{1,\ldots,n\}\), with \(|S|=k\), and let
\(x^S\in\{0,1\}^n\) denote its incidence vector. The corresponding
rank-one matrix $Z^S:=x^S(x^S)^\top$ is feasible for~\eqref{eq:fr_sdp}. Indeed, \(\diag(Z^S)=x^S\) and \(Z^S\one_n=x^S(x^S)^\top\one_n=kx^S=k\diag(Z^S)\). Moreover, \(\tr(Z^S)=\|x^S\|^2=k\), while \(Z^S\succeq 0\) and \(Z^S\geq\zero_{n\times n}\).
For a prescribed subset \(S\), we say that the relaxation is
\emph{tight at \(S\)} if \(Z^S\) is an optimal solution of
\eqref{eq:fr_sdp}. We say that \emph{exact recovery} occurs if
\(Z^S\) is the unique optimal solution. Note that, without uniqueness, the SDP may admit additional
fractional or integral optimizers and therefore may not identify the
subset \(S\) unambiguously.
Inspired by the primal--dual approach used in exact-recovery
analyses of SDP relaxations for clustering problems
\cite{iguchi2017probably,li2020birds}, we exploit the duality
result established in the previous section to derive a
necessary-and-sufficient dual characterization of optimality of the
candidate matrix \(Z^S\). We then strengthen this characterization to
obtain sufficient conditions for uniqueness. The dual certificate is
subsequently translated into explicit geometric requirements on the
data configuration. In particular, the resulting conditions separate
the role of the internal geometry of \(S\), measured through its
dispersion and anisotropy around the centroid, from the role of its
separation from its complement. Finally, by using the
necessary part of the primal--dual optimality conditions, we derive a
geometric condition that every tight configuration must satisfy.
Accordingly, the central question of this section is:
\begin{center}
\emph{Under what geometric conditions is \(Z^S\) the unique solution
of the SDP relaxation~\eqref{eq:fr_sdp}?}
\end{center}

The following lemma specializes the complementary slackness conditions of Lemma~\ref{lem:complementarity} to the candidate solution \(Z^S\).

\begin{lemma}
\label{lem:specialized_complementarity}
Let \((u,\tau,\Gamma,Q)\) be dual feasible. Then \(Z^S\) and \((u,\tau,\Gamma,Q)\) are primal--dual optimal if and only if $\Gamma_{SS}=\zero_{k\times k}$ and
$Qx^S=\zero_n$.
\end{lemma}

\begin{proof}
By Lemma~\ref{lem:complementarity}, \(Z^S\) and \((u,\tau,\Gamma,Q)\) are primal--dual optimal if and only if \(\Gamma\circ Z^S=\zero_{n\times n}\) and \(QZ^S=Z^SQ=\zero_{n\times n}\). Since \(Z^S_{ij}=1\) if \(i,j\in S\) and \(Z^S_{ij}=0\) otherwise, the condition \(\Gamma\circ Z^S=\zero_{n\times n}\) is equivalent to \(\Gamma_{SS}=\zero_{k\times k}\).
Moreover, $QZ^S=Qx^S(x^S)^\top$.
Hence, \(Qx^S=\zero_n\) implies \(QZ^S=\zero_{n\times n}\), and, since \(Q\) and \(Z^S\) are symmetric, also \(Z^SQ=(QZ^S)^\top=\zero_{n\times n}\). Conversely, if \(QZ^S=\zero_{n\times n}\), then $\zero_n
=
QZ^Sx^S
=
Qx^S(x^S)^\top x^S
=
kQx^S$,
where we used \((x^S)^\top x^S=k\). Therefore, \(Qx^S=\zero_n\), which proves the result.
\qed
\end{proof}

The next proposition gives sufficient conditions, expressed in terms of a single residual matrix, that guarantee not only optimality but also uniqueness of \(Z^S\).

\begin{proposition}\label{prop:tightness}
Suppose there exist \(u\in\mathbb R^n\) and \(\tau\in\mathbb R\) such that the matrix
\[
H:=
{\frac{1}{2k^2}D}
-
\operatorname{Diag}(u)
+
\frac{1}{2k}(u\mathbf{1}_n^\top+\mathbf{1}_nu^\top)
-
\tau \mathbf{I}_n
\]
satisfies
\begin{equation}\label{eq:block_psd}
H_{SS}\succeq0,
\end{equation}
\begin{equation}\label{eq:block_kernel}
H_{SS}\mathbf{1}_k=\mathbf{0}_k,
\end{equation}
and
\begin{equation}\label{eq:block_positive}
H_{ij}>0
\qquad
\forall (i,j)\notin S\times S.
\end{equation}
Then \(Z^S\) is the unique optimal solution of the SDP relaxation~\eqref{eq:fr_sdp}.
\end{proposition}

\begin{proof}
We first construct a dual feasible solution. Define matrices $\Gamma, Q \in \mathbb{S}^n$ by
\[
\Gamma_{ij}:=
\begin{cases}
0, & (i,j)\in S\times S,\\
H_{ij}, & (i,j)\notin S\times S
\end{cases}
\qquad\text{and}\qquad
Q:=H-\Gamma.
\]
The matrix \(\Gamma\) collects the entries of \(H\) outside the block \(S\times S\), while
\(Q\) retains only the principal block indexed by \(S\).
By~\eqref{eq:block_positive}, all entries of \(\Gamma\) are nonnegative, so $\Gamma\ge\mathbf{0}_{n \times n}$. Moreover, \(Q\) vanishes outside \(S\times S\) and 
$Q_{SS} = H_{SS}$, so~\eqref{eq:block_psd} implies $Q\succeq0$. 
The definition of \(H\) yields
\begin{equation*}
    \operatorname{Diag}(u)
-\frac{1}{2k}(u\mathbf{1}_n^\top+\mathbf{1}_nu^\top)
+\tau \mathbf{I}_n
+\Gamma
+Q
=
\frac{1}{2k^2}D,
\end{equation*}
and hence the point $(u, \tau, \Gamma, Q)$ is dual feasible. 
Furthermore, from the definition of $\Gamma$, we have $\Gamma_{SS} = \zero_{k \times k}$, while \eqref{eq:block_kernel} gives $Qx^S = \zero_n$. Therefore, by Lemma~\ref{lem:specialized_complementarity}, $Z^S$ and $(u, \tau, \Gamma, Q)$ are primal--dual optimal.

It remains to prove uniqueness. Let \(\bar Z \) be any optimal solution of the SDP. Since $(u, \tau, \Gamma, Q)$ is dual optimal, complementary slackness gives $\Gamma \circ \bar{Z} = \zero_{n \times n}$. By~\eqref{eq:block_positive}, $\Gamma_{ij} > 0$ for every $(i, j) \notin S \times S$. Hence, $\bar{Z}_{ij} = 0$ for all $(i, j) \notin S \times S$, so $\bar{Z}$ is supported entirely on \(S\times S\).
The affine constraints of~\eqref{eq:fr_sdp} then reduce to $\bar{Z}_{SS} \mathbf{1}_k=k \operatorname{diag}(\bar{Z}_{SS})$ and {$\operatorname{tr}(\bar Z_{SS})=k$} with \(\bar{Z}_{SS} \succeq 0\). Premultiplying the row-sum constraints by \(\mathbf{1}_k^\top\) gives $\mathbf{1}_k^\top \bar{Z}_{SS}\mathbf{1}_k
=
k\,\operatorname{tr}(\bar{Z}_{SS})
=
k^2$.
On the other hand, since \(\bar{Z}_{SS}\succeq0\), we have
\[
\mathbf{1}_k^\top \bar{Z}_{SS}\mathbf{1}_k
\leq
\lambda_{\max}(\bar{Z}_{SS})\|\mathbf{1}_k\|_2^2
\leq
\operatorname{tr}(\bar{Z}_{SS})\,k
=
k^2.
\]
Thus equality holds throughout. Hence \(\lambda_{\max}(\bar{Z}_{SS})=\operatorname{tr}(\bar{Z}_{SS})=k\), \(\mathbf{1}_k\) is an eigenvector associated with \(\lambda_{\max}(\bar{Z}_{SS})\), and all remaining eigenvalues of \(\bar{Z}_{SS}\) are zero. Therefore, $\bar{Z}_{SS}=\mathbf{1}_k\mathbf{1}_k^\top$. It follows that $\bar Z=x^S(x^S)^\top=Z^S$.
Hence \(Z^S\) is the unique optimal solution of~\eqref{eq:fr_sdp}.
\qed
\end{proof}

Note that, strict positivity in condition~\eqref{eq:block_positive} is required only to establish uniqueness. If instead $H_{ij}\geq 0$ for all $(i,j)\notin S\times S$, the same construction provides a dual certificate proving that $Z^S$ is optimal, but not necessarily uniquely optimal.
Proposition~\ref{prop:tightness} reduces the construction of a dual certificate to three conditions on the residual matrix \(H\). 
Since \(H\) is an affine correction of the distance matrix, these algebraic conditions can subsequently be translated into geometric conditions on the configuration of the data points. 

\subsection{From the dual certificate to geometry}
\label{sec:deterministic-tightness}
We now derive explicit geometric conditions under which the SDP relaxation is guaranteed to be tight. {Let $\mu_S$ denote} the centroid of the subset $S$ and let $\bar{S} = \{1, \dots, n\} \setminus S$ be the complement of $S$. Let
\[
    \mu_S:=\frac1k\sum_{i\in S}p_i, \qquad r_i := p_i - \mu_S, \qquad q_i := \|r_i\|^2, \qquad \bar{q} := \frac{1}{k}\sum_{i \in S} q_i.
\]
Thus, $q_i$ is the squared radius of point $p_i$ from the centroid of $S$, while $\bar q$ is the average squared radius of the selected points.
Let $R_S \in \mathbb{R}^{d \times k}$ be the matrix with columns $r_i$, in any fixed ordering of $S$, and let $q \in \mathbb{R}^k$ collect the values $q_i$ for all {$i\in S$}. Centering yields
two identities that will be used throughout the paper.

\begin{lemma}\label{lemma:centering_identities}
For every $i \in S$,
\begin{align}
\label{eq:centering_Dsum}
    \sum_{j \in S} D_{ij} = k(q_i + \bar{q}).
\end{align}
Moreover, for every $v \in \one_k^\perp$,
\begin{align}\label{eq:centering_DSS}
    v^\top D_{SS}v = -2 \|R_S v\|^2.
\end{align}
\end{lemma}
\begin{proof}
See Appendix \ref{appendix:centering_identities}.
\end{proof}

\noindent
To state our first geometric tightness result, we define the quantities
\[
\Delta_S:=\max_{i,j\in S}D_{ij},
\qquad
\Delta_\times:=\min_{i\in S,\ j\in\bar S}D_{ij},
\qquad
\Delta_{\bar S}:=\min_{\substack{i,j\in\bar S\\ i\neq j}}D_{ij}.
\]
The quantity $\Delta_S$ measures the squared diameter of the planted subset, while $\Delta_\times$ and $\Delta_{\bar S}$ quantify the minimum separation between $S$ and its complement, and within the complement itself.
When $|\bar S|=1$, we adopt the convention $\Delta_{\bar S}=+\infty$.

\begin{theorem}
\label{thm:spherical-tightness}
Assume that all points in $S$ are equidistant from the centroid $\mu_S$, i.e.,
$\|p_i-\mu_S\|^2 = \bar q$ for all $i \in S$. Suppose further that $\Delta_\times > \Delta_S$ and $\Delta_{\bar{S}} > \Delta_{S}$. Then the SDP relaxation~\eqref{eq:fr_sdp} is tight at $S$ and $Z^S$ is the unique optimal solution.
\end{theorem}

\begin{proof}
See Appendix~\ref{appendix:proof_spherical}.
\end{proof}

\begin{figure}[!ht]
    \centering
    \includegraphics[width=0.78\textwidth]{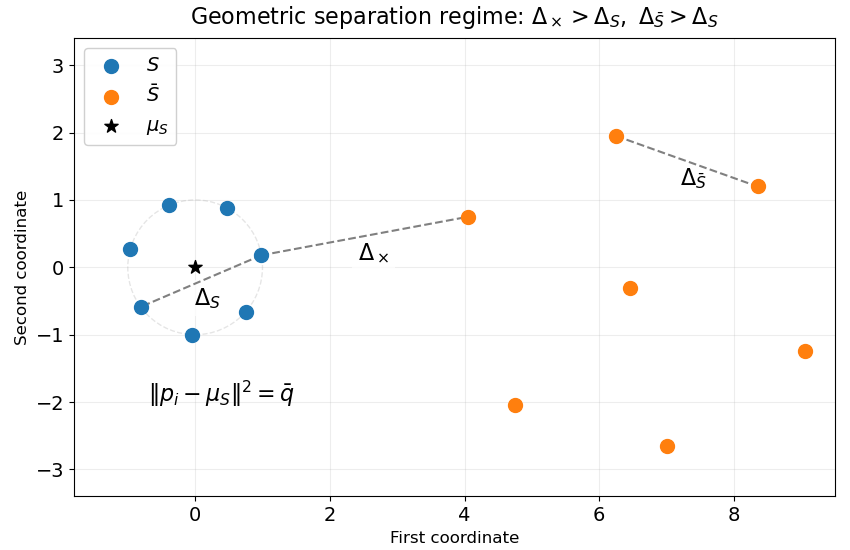}
    \caption{Schematic sufficient regime in Theorem~\ref{thm:spherical-tightness}: the selected points have a common radius and are separated from all unselected points.}
    \label{fig:tightness_delta}
\end{figure}

Theorem~\ref{thm:spherical-tightness} gives a simple and interpretable geometric condition for
exact recovery when the selected points are equidistant from their
centroid; see Figure~\ref{fig:tightness_delta}.

\begin{figure}[!ht]
    \centering
    \includegraphics[width=0.80\textwidth]{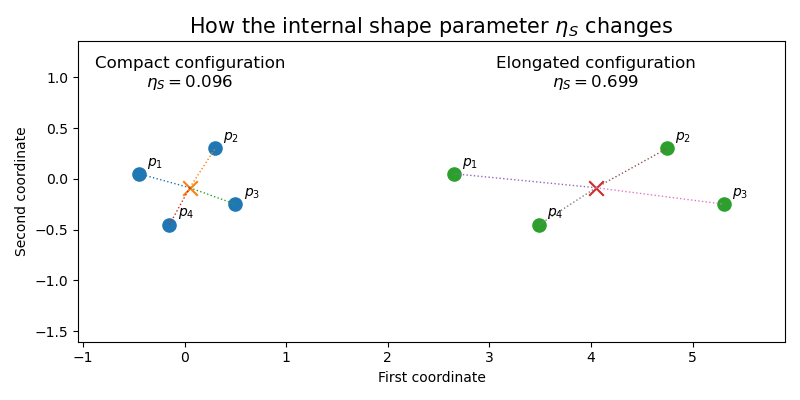}
    \caption{Illustration of the internal shape parameter: a compact configuration has a smaller value of $\eta_S$, whereas an elongated configuration has a larger value.}
    \label{fig:eta}
\end{figure}

We now remove this assumption by introducing an internal shape
parameter that accounts for the geometry of an arbitrary candidate subset.
Define
\begin{equation}
\eta_S
:=
\max_{\substack{v\in\mathbb{R}^k\\ v^\top \one_k=0,\ \|v\|=1}}
\left\{
\frac{1}{k}\sum_{i\in S}q_i v_i^2
+
\frac{1}{k^2}
\left\|
\sum_{i\in S}v_i r_i
\right\|^2
\right\}.
\label{eq:etaS}
\end{equation}
The quantity $\eta_S$ measures the internal geometric difficulty of certifying $S$ as optimal. It captures two distinct effects. The first term penalizes directions $v$ that place substantial weight on selected points lying far from the centroid. It is therefore sensitive to large individual radii and to nonuniform radial dispersion within $S$. The second term measures the extent to which the centered points combine coherently along a common direction, and therefore reflects anisotropy or elongation of the configuration. Consequently, $\eta_S$ is small when the selected points form a compact and approximately isotropic configuration, and becomes larger when the set is more elongated or contains points lying substantially farther from the centroid than the others. Figure~\ref{fig:eta} illustrates this behavior.

To compute \(\eta_S\), define
\[
B_S:=\frac{1}{k}\Diag(q)+\frac{1}{k^2}R_S^\top R_S.
\]
Let \(P:=\I_k-\frac{1}{k}\J_k\) be the orthogonal projector onto \(\mathbf{1}_k^\perp\). Let \(N\in\mathbb{R}^{k\times(k-1)}\) have orthonormal columns spanning \(\mathbf{1}_k^\perp\),  i.e., $N^\top N = \I_{k-1}$ and $N^\top \one_k = \zero_{k-1}$. Then \(NN^\top=P\), and every \(v\in\mathbf{1}_k^\perp\) can be written as \(v=Ny\) with \(\|v\|=\|y\|\). Hence
\[
\eta_S
= \max_{\substack{v^\top\mathbf{1}_k=0\\ \|v\|=1}}
v^\top B_Sv
=\max_{\norm{y}=1}y^\top N^\top B_SNy
=\lambda_{\max}(N^\top B_SN).
\]
Moreover, $PB_SP
=
N(N^\top B_SN)N^\top$.
Since \(N^\top N=\I_{k-1}\) and \(N^\top B_SN\succeq {0}\), the Rayleigh quotient gives $\lambda_{\max}(PB_SP)
=
\lambda_{\max}(N^\top B_SN)$.
Therefore,
\begin{equation}
\label{eq:eta_projected_eigenvalue}
\eta_S=\lambda_{\max}(PB_SP).
\end{equation}

To extend the tightness result beyond the equidistance assumption, define
the adjusted cross-separation
\begin{equation*}
\bar{\Delta}_\times
:=
\min_{\substack{i\in S, \ j\notin S}}
\left\{D_{ij}+q_i-\bar q\right\}.
\end{equation*}

The following result yields a sufficient condition for exact recovery that is computable from a candidate subset.

\begin{theorem}
\label{thm:checkable-certificate}
If
\begin{equation}
\label{eq:checkable_condition_unscaled}
\bar\Delta_\times>
k\left(\eta_S+\frac{\bar q}{k-1}\right)
\qquad\text{and}\qquad
\Delta_{\bar S}>
\frac{2k\bar q}{k-1},
\end{equation}
then \(Z^S\) is the unique optimal solution of the SDP relaxation~\eqref{eq:fr_sdp}. Consequently, \(S\) is the unique optimal solution of the MV$k$SP.
\end{theorem}

\begin{proof}
Condition~\eqref{eq:checkable_condition_unscaled} gives
\[
\frac{\bar q}{k-1}
<
\min\left\{
\frac{\Delta_{\bar S}}{2k},
\frac{\bar\Delta_\times}{k}-\eta_S
\right\}.
\]
Hence, we may choose a scalar $\theta$ such that
\begin{align}
\label{eq:small_s_def}
\frac{\bar q}{k-1}
<
\theta
<
\min\left\{
\frac{\Delta_{\bar S}}{2k},
\frac{\bar\Delta_\times}{k}-\eta_S
\right\}.
\end{align}
Define
\[
\tau := \frac{\bar{q}}{k}, \qquad u_i = -\eta_S + \frac{q_i - \bar q}{k} \quad (i \in S),  \qquad u_j = -\theta \quad (j \in \bar{S}).
\]
We now verify the three conditions in Proposition~\ref{prop:tightness}.
Since $\sum_{j \in S} u_j = -k \eta_S$, using \eqref{eq:centering_Dsum} from Lemma~\ref{lemma:centering_identities} gives, for every $i \in S$,
\begin{align}
    (H_{SS}\one_k)_i &= \frac{1}{2k^2}\sum_{j\in S} D_{ij} - \frac{1}{2} u_i + \frac{1}{2k}\sum_{j \in S} u_j - \tau \\
    &= \frac{q_i + \bar{q}}{2k} - \frac{1}{2}\left(-\eta_S + \frac{q_i - \bar q}{k} \right) - \frac{\eta_S}{2} - \frac{\bar q}{k} = 0.
\end{align}
Thus $H_{SS}\one_k = \zero_k$. 
To prove $H_{SS} \succeq 0$, let $v\in\mathbb{R}^k$ satisfy $v^\top\mathbf{1}_k=0$. The rank-one terms in the definition of $H$ vanish on such vectors. Hence
\begin{equation}
v^\top H_{SS}v
=
\frac{1}{2k^2}v^\top D_{SS}v
-\sum_{i\in S}u_i v_i^2
-\tau\|v\|^2.
\label{eq:vHv-start}
\end{equation}
Using \eqref{eq:centering_DSS} from Lemma~\ref{lemma:centering_identities} yields
\begin{equation}
\begin{aligned}
v^\top H_{SS}v
&=
-\frac{1}{k^2}\|R_{S} v\|^2
+
\eta_S\|v\|^2
-
\frac{1}{k}\sum_{i\in S}(q_i-\bar q)v_i^2
-
\frac{\bar q}{k}\|v\|^2\\
&=
\eta_S \| v\|^2 - v^\top B_S v \geq 0.
\end{aligned}
\label{eq:vHv-final}
\end{equation}
The last inequality is the definition of $\eta_S$, extended by homogeneity from unit vectors. Every vector 
\(y\in\mathbb R^k\) can be decomposed as $y=a\mathbf{1}_k+v$ with $v^\top \mathbf{1}_k = 0$, where $a = \frac{1}{k}\sum_{i=1}^k y_i$.
Since $H_{SS}\one_k = \zero_k$, this gives $y^\top H_{SS}y
=
v^\top H_{SS}v
\ge0$. Hence \(H_{SS}\succeq0\) holds.

It remains to verify that $H_{ij}>0$ for all $(i,j)\notin S\times S$.
First consider $i\in S$ and $j\notin S$. Since $i\neq j$, the diagonal terms vanish and
\[
H_{ij}
=
\frac{D_{ij} + q_i - \bar q - k(\eta_S + \theta)}{2k^2}.
\]
The numerator is at least $\bar{\Delta}_\times - k(\eta_S + \theta)$, which is positive by \eqref{eq:small_s_def}.
Next consider $i,j\notin S$ with $i\neq j$. Then
\begin{equation*}
H_{ij}
=
\frac{D_{ij} - 2k\theta}{2k^2} \geq  \frac{\Delta_{\bar S} - 2k\theta}{2k^2} > 0.
\end{equation*}
Finally, consider a diagonal entry $j\notin S$. Since $D_{jj}=0$, we obtain
\[
H_{jj}
=
-u_j+\frac{u_j}{k}-\tau
=
\frac{(k-1)\theta-\bar q}{k}.
\]
The lower bound on $\theta$ in \eqref{eq:small_s_def} gives
$(k-1)\theta>\bar q$, and hence $H_{jj}>0$.
Therefore $Z^S$ is the unique SDP optimizer, and $S$ is the unique MV$k$SP optimum.
\qed
\end{proof}

The certificate in Theorem~\ref{thm:checkable-certificate} separates the effects of internal geometry and external separation. Condition~\eqref{eq:checkable_condition_unscaled} requires \(\bar\Delta_\times\) to dominate the combined internal term \(k\bigl(\eta_S+\bar q/(k-1)\bigr)\), while \(\Delta_{\bar S}\) must exceed \(2k\bar q/(k-1)\). 
Thus, greater \(\eta_S\) requires stronger separation between selected and unselected points, while greater average dispersion within \(S\) strengthens both separation requirements.
Given a candidate set $S$, if Theorem~\ref{thm:checkable-certificate}  holds, then $S$ is certified as the unique MV$k$SP optimum without solving the combinatorial problem. Indeed, once the pairwise distances are available, the two separation quantities are obtained by scanning the relevant entries of $D$. The only spectral computation is the largest eigenvalue of
the $k\times k$ matrix $P B_S P$. 
Algorithm~\ref{alg:exact-recovery} summarizes the exact-recovery certificate provided by Theorem~\ref{thm:checkable-certificate}.

\begin{algorithm}[!ht]
\caption{Polynomial-time exact-recovery certificate}
\label{alg:exact-recovery}
\centering
\small
\renewcommand{\arraystretch}{1.15}
\begin{tabular}{@{}p{0.04\linewidth}p{0.89\linewidth}@{}}
1 &
Compute
$\mu_S=\frac{1}{k}\sum_{i\in S}p_i$,
$r_i=p_i-\mu_S$,
$q_i=\|r_i\|^2$, and
$\bar q=\frac{1}{k}\sum_{i\in S}q_i$.
\\[2pt]

2 &
Compute $\eta_S$ using
\eqref{eq:eta_projected_eigenvalue}.
\\[2pt]

3 &
Compute
$\bar{\Delta}_{\times}
 =\min_{i\in S,\;j\notin S}\{D_{ij}+q_i-\bar q\}$
and
$\Delta_{\bar S}
 =\min_{i,j\notin S,\;i\neq j}D_{ij}$.
\\[2pt]

4 &
If
$\bar\Delta_\times>
k\left(\eta_S+\frac{\bar q}{k-1}\right)$
and
$\Delta_{\bar S}>
\frac{2k\bar q}{k-1}$,
certify that $S$ is the unique optimum. Otherwise, the test is inconclusive.
\\
\end{tabular}
\end{algorithm}

We now derive a necessary condition for tightness. The proof has two
ingredients. First, positive semidefiniteness of the active dual block
forces a scalar dual parameter to dominate the internal shape parameter
\(\eta_S\). Second, the dual equations associated with each outside point
convert this scalar bound into a lower bound on its distance from the
centroid of \(S\).

\begin{proposition}
\label{lem:active-rho}
Suppose that \(Z^S\) is optimal for the SDP relaxation, and let
\((u,\tau,\Gamma,Q)\) be any optimal dual solution. Define
$\rho:=\tau-\frac1k \one_k^\top u_S$.
Then
\[
\tau=\frac{\bar q}{k},
\qquad
u_S=\frac{q}{k}-\rho\one_k,
\qquad
Q_{SS}
=
\rho P
-
P\left(
B_S+\frac{\bar q}{k}I_k
\right)P.
\]
Consequently, $\rho\geq\eta_S+\frac{\bar q}{k}$.
\end{proposition}

\begin{proof}
Complementary slackness in Lemma~\ref{lem:specialized_complementarity} gives $\Gamma_{SS}=\mathbf{0}_{k \times k}$ and $Q_{SS}\mathbf{1}_k=\mathbf{0}_k$. Let $C = D / (2k^2)$.
The $S \times S$ block of the dual equality is therefore

\begin{equation} \label{eq:active_dual_block}
Q_{SS}
=
C_{SS}
-\Diag(u_S)
+\frac1{2k}
\left(
u_S\mathbf{1}_k^\top+\mathbf{1}_ku_S^\top
\right)
-\tau \mathbf{I}_k.
\end{equation}
Post-multiplying \eqref{eq:active_dual_block} by $\mathbf{1}_k$ and using \(Q_{SS}\mathbf{1}_k=\mathbf{0}_k\) gives
\begin{equation}\label{eq:active_dual_block_sum}
\zero_k
=
C_{SS}\one_k
-\frac12u_S
+\frac{\one_k^\top u_S}{2k}\one_k
-\tau \one_k.
\end{equation}
Pre-multiplying \eqref{eq:active_dual_block_sum} by $\one_k^\top$ and using the centroid identity $\one_k^\top D_{SS} \one_k
=
2k\sum_{i\in S}\|p_i-\mu_S\|^2
=
2k^2\bar q$ yields $0 = \bar q - k \tau$, hence $\tau=\frac{\bar q}{k}$. Next, using the definition of $\rho$, equation \eqref{eq:active_dual_block_sum} yields $u_S = 2 C_{SS} \one_k - (\tau + \rho) \one_k = \frac{1}{k^2} D_{SS} \one_k - (\tau + \rho) \one_k$. By Lemma \ref{lemma:centering_identities}, $D_{SS}\one_k = k(q + \bar{q} \one_k)$, so this simplifies to $u_S = \frac{q}{k} - \rho \one_k$.

Since $Q_{SS}\one_k = \zero_k$ and $Q_{SS}$ is symmetric, we have $Q_{SS} = PQ_{SS}P$. Multiplying \eqref{eq:active_dual_block} on both sides by $P$ and using $P \one_k = \zero_k$ and $P^2 = P$ therefore gives 
\begin{equation}\label{eq:projected_Q}
    Q_{SS} = PC_{SS}P - P \Diag(u_S)P - \tau P.
\end{equation}
The squared-distance identity $D_{SS} = q\one_k^\top + \one_k q^\top - 2 R_S^\top R_S$ and $R_S \one_k = \zero_d$ imply $PC_{SS}P = -\frac{1}{k^2}R_S^\top R_S$ and $P(R_S^\top R_S)P = R_S^\top R_S$. Substituting the expression for $u_S$ and $\tau$ into \eqref{eq:projected_Q}, we obtain
\[
Q_{SS} = -\frac{1}{k^2} R_S^\top R_S - \frac{1}{k} P \Diag(q)P + \rho P - \frac{\bar{q}}{k}P = \rho P - P \left(B_S + \frac{\bar{q}}{k} \mathbf{I}_k \right) P.
\]

Finally, \(Q\succeq0\) implies $Q_{SS}\succeq0$.
Let \(v\in\mathbb R^k\) satisfy $v^\top \one_k=0$ and $\|v\| = 1$.
Since \(v\in \one_k^\perp\), the projector \(P\) satisfies $Pv=v$. Thus,
\[
0 \leq v^\top Q_{SS}v
=
\rho v^\top Pv
-
v^\top
P\left(
B_S+\frac{\bar q}{k}\I_k
\right)Pv =
\rho
-v^\top B_Sv
-\frac{\bar q}{k}.
\]
Thus, $\rho
\geq
v^\top B_Sv+\frac{\bar q}{k}$
for every unit vector \(v\perp \one_k\). Taking the maximum of $v^\top B_S v$ over all such vectors gives $\rho\geq\eta_S+\frac{\bar q}{k}$, completing the proof.
\qed
\end{proof}

We now translate the conditions in Proposition~\ref{lem:active-rho} into geometric conditions that are necessary for tightness.

\begin{theorem}
\label{thm:necessary-tightness}
Suppose that the SDP relaxation is tight at $S$, i.e., that $Z^S$ is optimal. Then every $\ell\notin S$ satisfies
\[
\|p_\ell-\mu_S\|^2
\geq
k\eta_S+\frac{\bar q}{k-1}.
\]
\end{theorem}

\begin{proof}
Fix $\ell\notin S$, let $(u,\tau,\Gamma,Q)$ be an optimal dual
solution, and set $C:=D/(2k^2)$.
Define $\delta_\ell:=\Gamma_{\ell\ell}+Q_{\ell\ell}\geq0$.
Since $C_{\ell\ell}=0$, the diagonal dual equality gives $0
=
\frac{k-1}{k}u_\ell+\tau+\delta_\ell$,
hence $u_\ell
=
-\frac{k}{k-1}(\tau+\delta_\ell)$.
For \(i\in S\), the $(\ell, i)$ entry of the dual equality is
\[
C_{\ell i}
=
-\frac{u_\ell+u_i}{2k}
+\Gamma_{\ell i}+Q_{\ell i}.
\]
Summing over $i\in S$, using
$\sum_{i\in S}Q_{\ell i}=0$ from complementary slackness
and $\frac{1}{k}\sum_{i\in S}u_i=\tau-\rho$, we obtain
\begin{align*}
2\sum_{i\in S}C_{\ell i}
&=-u_\ell-\frac{1}{k}\sum_{i\in S}u_i
  +2\sum_{i\in S}\Gamma_{\ell i}\\
&=\rho+\frac{\tau}{k-1}
  +\frac{k}{k-1}\delta_\ell
  +2\sum_{i\in S}\Gamma_{\ell i}\\
&\geq \rho+\frac{\tau}{k-1}.
\end{align*}
Proposition~\ref{lem:active-rho} gives
$\rho\geq\eta_S+\bar q/k$ and $\tau=\bar q/k$, so
\begin{equation}\label{eq:sum_bound}
2\sum_{i\in S}C_{\ell i}
\geq\eta_S+\frac{\bar q}{k-1}.
\end{equation}
On the other hand, the centroid identity gives
\begin{equation}\label{eq:centr_iden}
2\sum_{i\in S}C_{\ell i}
=\frac{1}{k^2}\sum_{i\in S}D_{\ell i}
=\frac{1}{k}\left(\|p_\ell-\mu_S\|^2+\bar q\right).
\end{equation}
Combining \eqref{eq:sum_bound} and \eqref{eq:centr_iden}
and rearranging yields $\|p_\ell-\mu_S\|^2
\geq k\eta_S+\frac{\bar q}{k-1}$.
\qed
\end{proof}

We remark that 
for $\ell\notin S$ and $i\in S$, let $\omega_{\ell i} := D_{\ell i} + q_i - \bar q$, so that $\bar\Delta_\times = \min_{\ell\notin S,\, i\in S}\omega_{\ell i}$. Since $\sum_{i\in S}(q_i-\bar q)=0$, identity~\eqref{eq:centr_iden} gives $\frac1k\sum_{i\in S}\omega_{\ell i} = \|p_\ell-\mu_S\|^2 + \bar q$. Theorem~\ref{thm:necessary-tightness} is therefore equivalent to
\[
\frac1k\sum_{i\in S}\omega_{\ell i} \;\ge\; k\Bigl(\eta_S + \frac{\bar q}{k-1}\Bigr)
\qquad \forall \ell\notin S,
\]
whereas the first  condition in~\eqref{eq:checkable_condition_unscaled} requires
\[
\min_{i\in S}\omega_{\ell i} \;>\; k\Bigl(\eta_S + \frac{\bar q}{k-1}\Bigr)
\qquad \forall \ell\notin S.
\]
The sufficient and necessary conditions thus share the same threshold. They differ only in that the former aggregates the separation from $S$ by a minimum rather than an average, and requires separation among the unselected points.

\section{Branch-and-cut algorithm}
\label{sec:algorithm}


In this section, we describe the key ingredients of the algorithm for solving problem \eqref{eq:quadratic_distance} to global optimality.
\subsection{Lower bound computation}
Doubly nonnegative (DNN) relaxations are known to significantly strengthen standard
SDPs for many combinatorial optimization problems, as they
simultaneously exploit positive semidefiniteness and entrywise nonnegativity of the
lifted matrix \cite{kim2022doubly}. However, solving DNNs with interior-point methods becomes computationally
prohibitive for large-scale instances due to their high memory and time complexity.
To address this limitation, first-order methods based on augmented Lagrangian
frameworks have been proposed, which are able to handle larger problem sizes by
trading off accuracy for scalability \cite{yang2015sdpnal,sun2020sdpnal+}.

The relaxation can be further tightened by incorporating valid linear inequalities
derived from the Boolean quadric polytope (BQP), i.e., the convex hull of all
points $(x, X)$ satisfying $x \in \{0,1\}^n$ and $X = xx^\top$. These inequalities
encode relationships between $x_i$ and the lifted variables $X_{ij}$, and are
known to strengthen SDP relaxations significantly
\cite{hrga2021madam,gusmeroli2022biqbin,locatelli2025fix}.
In the facially reduced formulation, $X=Z$ and
$x=\frac{1}{k}Z\mathbf{1}_n=\operatorname{diag}(Z)$.
Thus, these inequalities can be expressed directly in terms of $Z$.
We consider the following well-known classes.

\paragraph{McCormick inequalities.}
For all distinct $i, j \in \{1, \dots, n\}$, the Reformulation--Linearization Technique (RLT) yields the
constraints
\begin{align}
Z_{ij} &\le Z_{ii}, \label{eq:rlt1}\\
Z_{ij} &\le Z_{jj}, \label{eq:rlt2}\\
Z_{ij} &\ge Z_{ii}+Z_{jj}-1. \label{eq:rlt3}
\end{align}

\paragraph{Triangle inequalities.}
Additional valid inequalities are obtained by considering triples of variables.
For any distinct indices $i,j,h \in \{1,\dots,n\}$, these read
\begin{align}
Z_{ij}+Z_{ih}-Z_{jh} &\le Z_{ii}, \label{eq:tri1}\\
Z_{ij}+Z_{jh}-Z_{ih} &\le Z_{jj}, \label{eq:tri2}\\
Z_{ih}+Z_{jh}-Z_{ij} &\le Z_{hh}, \label{eq:tri3}\\
Z_{ij}+Z_{ih}+Z_{jh} &\ge Z_{ii}+Z_{jj}+Z_{hh}-1. \label{eq:tri4}
\end{align}

Incorporating all inequalities
\eqref{eq:rlt1}--\eqref{eq:rlt3} and
\eqref{eq:tri1}--\eqref{eq:tri4} explicitly is computationally
prohibitive due to their large number; in particular, there are
$\mathcal{O}(n^3)$ triangle inequalities.
Instead, we add them in a cutting-plane fashion: starting from the basic
SDP relaxation, we iteratively identify violated inequalities through
a separation procedure and add them dynamically to the relaxation.
In order to scale to large instances, we solve the SDP relaxation
\eqref{eq:fr_sdp} strengthened with inequalities  
\eqref{eq:rlt1}--\eqref{eq:tri4} by a first-order method. Since the returned primal--dual
iterate is generally only approximately feasible, the dual objective value
cannot be used directly as a certified lower bound. Therefore, we apply the LP-based post-processing step of \cite{cerulli2021improving}
that recovers a dual feasible solution and hence a valid lower bound. The strengthened relaxation can be written compactly as
\begin{equation}
\label{eq:fr_sdp_compact_ineq}
\begin{aligned}
\min_{Z \in \mathbb{S}^n}\quad & \frac{1}{2k^2}\langle D, Z\rangle \\
\text{s.t.}\quad
& \mathcal{A}(Z) = \mathbf{0}_n, \ \mathcal{B}(Z) \geq b,\ \langle \mathbf{I}_n, Z\rangle = k,\\
& Z \succeq 0,\ Z \geq \mathbf{0}_{n \times n},
\end{aligned}
\end{equation}
where the linear operator \(\mathcal{A}:\mathbb{S}^n \to \mathbb{R}^n\) is
defined by $\mathcal{A}(Z) = \diag(Z) - \frac{1}{k} Z\mathbf{1}_n$.
The adjoint operator \(\mathcal{A}^*:\mathbb{R}^n \to \mathbb{S}^n\) is given by $\mathcal{A}^*(u)
= \Diag(u) - \frac{1}{2k}\left(u\mathbf{1}_n^\top + \mathbf{1}_n u^\top\right)$ for some $u \in \mathbb{R}^n$. Finally, the linear operator \(\mathcal{B}:\mathbb{S}^n \to \mathbb{R}^l\)
encodes RLT inequalities
\eqref{eq:rlt1}--\eqref{eq:rlt3} and triangle inequalities \eqref{eq:tri1}--\eqref{eq:tri4}, with adjoint operator
\(\mathcal{B}^*:\mathbb{R}^l \to \mathbb{S}^n\) and right-hand-side vector $b \in \mathbb{R}^l$.

Let \(u \in \mathbb{R}^n\) be the multiplier associated with
\(\mathcal{A}(Z)=\mathbf 0_n\), let \(\pi\in\mathbb R^l_+\) be the multiplier
associated with \(\mathcal{B}(Z)\ge b\), let \(\tau\in\mathbb R\) be the
multiplier associated with \(\langle \mathbf{I}_n,Z\rangle=k\), let
\(\Gamma\in\mathbb S^n\) be the multiplier associated with \(Z\ge\mathbf{0}_{n \times n}\), and let
\(Q\in\mathbb S^n\) be the multiplier associated with \(Z\succeq0\). The dual
of problem \eqref{eq:fr_sdp_compact_ineq} is
\begin{equation}
\label{eq:fr_compact_dual_ineq}
\begin{aligned}
\max_{u,\pi,\tau,\Gamma,Q}\quad
& k \tau + b^\top \pi \\
\text{s.t.}\quad
&
\mathcal{A}^*(u)
+
\mathcal{B}^*(\pi)
+
\tau \mathbf{I}_n
+
\Gamma
+
Q
=
\frac{1}{2 k^2} D, \\
& Q\succeq0, \ \pi\ge\mathbf{0}_l,\ \Gamma \ge\mathbf{0}_{n \times n}.
\end{aligned}
\end{equation}

By weak duality, every feasible solution of problem \eqref{eq:fr_compact_dual_ineq} provides a
valid lower bound for problem \eqref{eq:fr_sdp_compact_ineq}. 
Let $(\hat u,\hat \pi,\hat\tau,\hat \Gamma,\hat Q)
$ denote the approximate dual solution returned by the first-order solver.
To recover a certified lower bound, we first project \(\hat Q\) onto the
positive semidefinite cone. Let $\hat Q = \bar{U}\Lambda \bar{U}^\top
$ be an eigendecomposition of \(\hat Q\) with eigenvalues $\lambda_1, \dots, \lambda_n$, and define
$\widetilde Q := \bar{U}\Lambda_+\bar{U}^\top$, where $
\Lambda_+ := \Diag\bigl(\max\{\lambda_i,0\}\bigr)$.
Then \(\widetilde Q \succeq 0\) by construction. Keeping \(\widetilde Q\) fixed, we solve the linear program
\begin{equation}
\label{eq:fr_post_lp}
\begin{aligned}
\max_{u,\pi,\tau,\Gamma}\quad
& k \tau + b^\top \pi \\
\text{s.t.}\quad
&
\mathcal{A}^*(u)
+
\mathcal{B}^*(\pi)
+
\tau \mathbf{I}_n
+
\Gamma
=
\frac{1}{2 k^2} D-\widetilde Q, \\
& \pi\ge\mathbf{0}_l,\ \Gamma\ge\mathbf{0}_{n \times n}.
\end{aligned}
\end{equation}

Let $(\widetilde u,\widetilde \pi,\widetilde\tau,\widetilde\Gamma)$
be an optimal solution of problem~\eqref{eq:fr_post_lp}.
Then $(\widetilde u,\widetilde \pi,\widetilde\tau,
\widetilde\Gamma,\widetilde Q)$ is feasible for
problem~\eqref{eq:fr_compact_dual_ineq}, so
$k\widetilde\tau+b^\top\widetilde \pi$ provides a valid lower bound
for the original integer problem.

\subsection{Branching and node subproblems}
Using cutting planes to tighten the SDP relaxation may not be sufficient to find or prove the optimality
of the best solution found so far. If there are no violated cuts or the bound does not
improve substantially after adding valid inequalities, we terminate the
cutting plane generation and proceed to branching.
We branch on a fractional diagonal entry \(Z_{ii}\), creating children with \(Z_{ii}=0\) and \(Z_{ii}=1\).
The next result shows that these branching decisions admit strong propagation and lead to a reduction of the problem dimension. 

\begin{lemma}\label{lemma:propagation}
Let $Z\in\mathbb{S}^n$ be feasible for \eqref{eq:fr_sdp}. Fix $i\in\{1,\dots,n\}$.
\begin{enumerate}
    \item If $Z_{ii}=0$, then $Z_{ij}=0$ for all $j\in\{1,\dots,n\}$.
    \item If $Z_{ii}=1$, then $Z_{ij}=Z_{jj}$ for all $j\in\{1,\dots,n\}$.
\end{enumerate}
Consequently, index $i$ can be eliminated in both cases, reducing the problem dimension from $n$ to $n-1$. In the case $Z_{ii}=0$, the cardinality remains $k$, while in the case $Z_{ii}=1$, it reduces to $k-1$.
\end{lemma}

\begin{proof}
See Appendix~\ref{appendix:lemma_propagation}.
\end{proof}

This recursive structure motivates the definition of a generic subproblem at each node of the
branch-and-bound tree. In our framework, a node is characterized by two disjoint index sets $F_1,F_0 \subseteq \{1,\dots,n\}$,
where $F_1$ contains indices fixed to $1$ and $F_0$ contains indices fixed to $0$. 
Let $U = \{1,\dots,n\}\setminus (F_1\cup F_0)$
denote the set of unfixed indices, and define $\bar{n} = |U|$ and $\bar{k} = k-|F_1|$.
Thus, at the current node, exactly $\bar{k}$ additional points must still be selected among the indices in $U$.
We introduce binary variables $y \in \{0,1\}^{\bar n}$ for the unfixed indices, so that $(x_{F_1},\,x_{F_0},\,x_U)=(\mathbf{1},\,\mathbf{0},\,{y})$.
The cardinality constraint becomes $\mathbf{1}^\top y = \bar{k}$. Define the vector $\delta_U \in \mathbb{R}^{\bar{n }}$ such that, for all $j \in U$, we have
\[
(\delta_U)_j = \sum_{i\in F_1} D_{ij},
\]
and the constant
\[
c(F_1)
=
\frac{1}{2k^2}
\sum_{i\in F_1}\sum_{j\in F_1} D_{ij}.
\]
The objective function can be expanded as
\begin{align*}
\frac{1}{2k^2}x^\top D x
&=
\frac{1}{2k^2}
\sum_{i,j\in F_1} D_{ij}
+
\frac{1}{k^2}
\sum_{i\in F_1}\sum_{j\in U} D_{ij} y_j
+
\frac{1}{2k^2}
\sum_{i,j\in U} D_{ij} y_i y_j \\
&=
c(F_1)
+
\frac{1}{k^2}\delta_U^\top y
+
\frac{1}{2k^2} y^\top D_{UU} y.
\end{align*}
Therefore, the exact subproblem at a node is
\begin{equation}
\label{eq:node_bqp}
\begin{aligned}
\min_{y \in \{0,1\}^{\bar n}} \quad
&
c(F_1)
+
\frac{1}{k^2}\delta_U^\top y
+
\frac{1}{2k^2} y^\top D_{UU} y \\
\text{s.t.}\quad
& \mathbf{1}_{\bar{n}}^\top y = \bar{k}.
\end{aligned}
\end{equation}
At the root node, $F_1 = F_0 = \emptyset$, which implies
$\bar{n} = n$, $\bar{k} = k$, $c(F_1)=0$, and $\delta_U=\mathbf{0}_{\bar{n}}$.
Hence, problem \eqref{eq:quadratic_distance} coincides with \eqref{eq:node_bqp}.
The same facial reduction argument presented in Section~\ref{sec:sdp} applies, replacing $k$ with $\bar{k}$. Indeed, the facially reduced SDP at the node reads as
\begin{equation}
\label{eq:node_fr_sdp}
\begin{aligned}
\min_{Z\in\mathbb{S}^{\bar n}} \quad
&
c(F_1)
+
\frac{1}{k^2}\delta_U^\top \operatorname{diag}(Z)
+
\frac{1}{2k^2}\langle D_{UU}, Z\rangle \\
\text{s.t.}\quad
& Z\mathbf{1}_{\bar n} = \bar{k} \operatorname{diag}(Z), \\
& \operatorname{tr}(Z) = \bar{k}, \\
& Z\succeq 0, \ Z \geq \mathbf{0}_{\bar{n} \times \bar{n}}.
\end{aligned}
\end{equation}
Every branching step therefore produces a smaller SDP of the same
form, with dimension equal to the number of unfixed indices.

\section{Computational results}
\label{sec:experiments}
In this section, we evaluate both the recovery theory and the computational
performance of the proposed branch-and-cut algorithm. Section~\ref{sub:exact} uses controlled synthetic
instances to assess the sufficient and necessary conditions for exact recovery.
Section~\ref{sub:bb_res} compares the exact algorithm with general-purpose optimization solvers on public benchmark datasets.

\subsection{Numerical assessment of exact recovery}
\label{sub:exact}
In this section we investigate the robustness of the tightness conditions on controlled synthetic instances. In particular, we study
how the sufficient and necessary conditions behave when an initially well-separated planted
subset that is uniquely optimal is progressively perturbed. Our experiments are designed to answer the following questions.

\begin{enumerate}
    \item How much perturbation can a certified instance tolerate before the
    sufficient condition ceases to hold?

    \item Does failure of the sufficient condition coincide with failure of
    SDP tightness, or is there a region in which the SDP remains
    exact although our certificate no longer applies?

    \item At what perturbation level does the planted subset cease to be the
    globally optimal solution of the original combinatorial problem?


\end{enumerate}

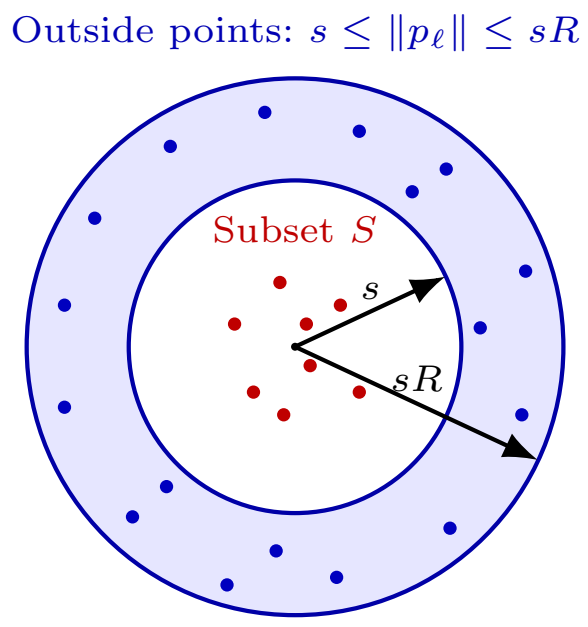
\begin{figure}[!ht]
\centering
\begin{tikzpicture}[
  scale=0.5,
  every node/.style={font=\scriptsize}
]

\def\innerR{2.2}
\def\outerR{3.55}

\fill[blue!10] (0,0) circle (\outerR);
\fill[white] (0,0) circle (\innerR);

\draw[blue!65!black, thick] (0,0) circle (\innerR);
\draw[blue!65!black, thick] (0,0) circle (\outerR);

\foreach \x/\y in {
 -0.80/0.30,
 -0.55/-0.60,
 -0.20/0.85,
  0.20/-0.25,
  0.60/0.55,
  0.85/-0.60,
 -0.15/-0.90,
  0.15/0.30
}{
  \fill[red!75!black] (\x,\y) circle (2.5pt);
}

\foreach \x/\y in {
 -3.05/0.55,
 -2.65/1.70,
 -1.65/2.65,
 -0.40/3.10,
  0.85/2.85,
  2.00/2.35,
  3.05/1.00,
  3.00/-0.90,
  2.05/-2.40,
  0.55/-3.05,
 -0.90/-3.15,
 -2.15/-2.25,
 -3.05/-0.80,
  1.55/2.05,
 -1.70/-1.85,
  2.45/0.25,
 -0.25/-2.70
}{
  \fill[blue!75!black] (\x,\y) circle (2.5pt);
}

\fill (0,0) circle (1.5pt);

\draw[-{Latex}, thick]
  (0,0) --
  ({\innerR*cos(25)},{\innerR*sin(25)})
  node[midway, above, inner sep=2pt] {\(s\)};

\draw[-{Latex}, thick]
  (0,0) --
  ({\outerR*cos(-25)},{\outerR*sin(-25)})
  node[midway, above, inner sep=2pt] {\(sR\)};

\node[red!75!black, align=center]
  at (0,1.55)
  {Subset \(S\)};

\node[blue!70!black, align=center]
  at (0,4.15)
  {Outside points: \(s\leq\|p_\ell\|\leq sR\)};

\end{tikzpicture}
\caption{Schematic two-dimensional representation of the hollow spherical
data model.}
\label{fig:hollow-ball-model}
\end{figure}

\paragraph{Data generation.}
Let $S \subset \{1, \dots, n\}$ denote the planted subset with $|S| = k$. We draw
\(z_i\sim N(0,I_d)\) independently, define
\(\bar z_S=\frac{1}{k}\sum_{j\in S}z_j\), and set
\[
 p_i=
 \frac{z_i-\bar z_S}
 {\sqrt{\frac{1}{k}\sum_{j\in S}\|z_j-\bar z_S\|^2}},
 \qquad i\in S.
\]
The planted centroid is therefore the origin and its average squared
radius is \(\bar q=\frac{1}{k}\sum_{i\in S}\|p_i\|^2=1\).  The planted points are not constrained to have equal norms, so the generator retains the random internal asymmetry measured by \(\eta_S\).
The remaining \(m=n-k\) points are sampled from a hollow spherical
region centered at the same origin. These points fill the volume between two concentric hyperspheres; in two dimensions, the corresponding region is a ring, as illustrated in
Figure~\ref{fig:hollow-ball-model}.

For each \(\ell\notin S\), we independently draw
\(g_\ell\sim N(0,I_d)\) and \(W_\ell\sim\operatorname{Unif}(0,1)\), and set \(e_\ell=g_\ell/\|g_\ell\|\) and $R_\ell=\left[1+W_\ell(R^d-1)\right]^{1/d}$.
The outside point is then \(p_\ell=sR_\ell e_\ell\), and hence \(s\leq\|p_\ell\|\leq sR\).  
Consequently, the outside points are uniformly distributed with respect to \(d\)-dimensional volume in the hollow spherical region. The parameter \(R>1\) controls the relative thickness of this region, while \(s\) sets its inner radius and overall radial scale.  As \(m\) varies, we choose \(R\) so that the volume per outside point of the unscaled region, whose inner and outer radii are \(1\) and \(R\), remains constant.  At fixed dimension \(d\), this volume is proportional to \(R^d-1\); therefore, we set
\(R=[1+(m/m_{\rm ref})(R_{\rm ref}^d-1)]^{1/d}\), where
\(m_{\rm ref}\) and \(R_{\rm ref}\) define the reference configuration. The common
scale \(s\) is selected separately as described below.

\paragraph{Calibration and reference-instance selection.}
The scale \(s\) determines the separation between the planted and outside points. Choosing it arbitrarily could make the zero-noise
problem either uncertified or trivially easy. We therefore calibrate \(s\) for every tested configuration using
\(N_{\mathrm{cal}}\) independent trials.
Each trial contains a newly sampled planted cloud and newly sampled outside directions
and relative radii. Within a trial these quantities remain fixed while $s$ varies, so changing $s$ only
moves the outside points radially.
For trial $j$ and scale $s$, we compute $
\kappa_j(s)
=
\min\left\{
\frac{\bar\Delta_\times}
     {k\left(\eta_S+\bar q/(k-1)\right)},
\frac{\Delta_{\bar S}}
     {2k\bar q/(k-1)}
\right\}$. The sufficient condition in Theorem~\ref{thm:checkable-certificate} holds if and only if \(\kappa_j(s)>1\). Starting from \(s=1\), we increase the scale until $\operatorname{median}_{1\leq j\leq N_{\mathrm{cal}}} \kappa_j(s) \geq \kappa_\star$,
where \(\kappa_\star>1\) is a prescribed target. Bisection is then used to
identify the smallest scale satisfying this inequality. 
After freezing \(s\), we generate independent candidate test instances. A candidate is retained only if the sufficient condition holds and the planted subset is the unique global optimum. Each retained dataset is a reference instance and it
serves as a starting point for the noise experiment. 

\paragraph{Noise experiment.}
Starting from each reference instance, we progressively increase an additive perturbation of the planted points while keeping the outside points fixed.
Let \(B\) be the number of reference instances and let \(T\) be the
number of independent perturbation patterns generated for each reference
instance. Furthermore, let \(P^{(b)}\in\mathbb R^{n\times d}\) denote reference instance
\(b\in\{1,\ldots,B\}\), with one observation in each row.  For each
\(t\in\{1,\ldots,T\}\), we independently draw a perturbation matrix
\(\Xi^{(b,t)}\in\mathbb R^{n\times d}\).  Its \(i\)-th row satisfies $ \xi_i^{(b,t)}\sim N(0,I_d/d)$ for all $ i\in S$
and $\xi_i^{(b,t)}=0$ for all $i\notin S$.
At noise level
\(\noise\geq0\), the perturbed dataset is $P^{(b,t)}(\noise)
 =
 P^{(b)}+\noise\Xi^{(b,t)}$.
For every planted point, $\mathbb E[\|\noise\xi_i^{(b,t)}\|^2]=\noise^2$. Therefore, \(\noise\) is the root-mean-square magnitude of its
additive perturbation, measured relative to an initial planted-cloud radius of one.
For fixed $b$ and $t$, the matrix $\Xi^{(b,t)}$ is drawn once
and held fixed as $\noise$ varies, defining a perturbation path.
Thus, increasing $\noise$ strengthens the same perturbation pattern,
while different values of $t$ correspond to independent patterns
applied to the same reference instance.
We evaluate $L$ equally spaced noise levels
$\sigma\in\{0,\sigma_{\max}/(L-1),\dots,\sigma_{\max}\}$.

\paragraph{Optimization and classification.}
At every perturbed instance, the original combinatorial problem is solved globally twice with $\tt{GUROBI}$.  The first solve is unrestricted.  The
second adds $\sum_{i\in S}x_i\leq k-1$,
which excludes exactly the planted subset and returns its best competitor; the resulting objective gap tests uniqueness.
We solve the SDP relaxation with \texttt{MOSEK}, using a default tolerance of $10^{-8}$ for primal and dual feasibility and the complementarity gap.
We classify the SDP as value-exact when its objective coincides
with the globally optimal combinatorial value to relative
tolerance $10^{-6}$. We declare recovery of the planted subset only when, in addition, the planted subset is uniquely optimal and the relative Frobenius distance between the SDP optimizer and the planted matrix \(x_Sx_S^\top\) is at most \(10^{-6}\).




\medskip
Throughout the reported experiments, we fix \(n=50\), \(d=10\), and \(k=5\). This setting allows for repeated solution of both the SDP relaxation and the original combinatorial problem in a reasonable amount of time. For the hollow-region construction, we use \(m_{\rm ref}=40\) and \(R_{\rm ref}=1.30\).  We use \(N_{\rm cal}=10\) calibration trials and target a median ratio \(\kappa_\star=1.5\). Since certification requires \(\kappa>1\), this provides a moderate initial margin while keeping the instances close enough to the boundary to exhibit informative transitions. We select the smallest \(s\) attaining this target, thereby avoiding excessive separation that would make the initial instances artificially easy. The perturbation experiment uses
\((B,T,L,\noise_{\max})=(10,10,16,3)\), corresponding to 10 reference instances, 10 independent perturbation paths per reference instance, and 16 noise levels between \(0\) and \(3\). Because all paths associated with the same reference instance coincide at \(\noise=0\), the unperturbed instance is solved only once. The experiment therefore contains \(B[1+T(L-1)]=1510\) distinct instances.

\begin{figure}[!ht]
    \centering
    \includegraphics[scale=0.54]{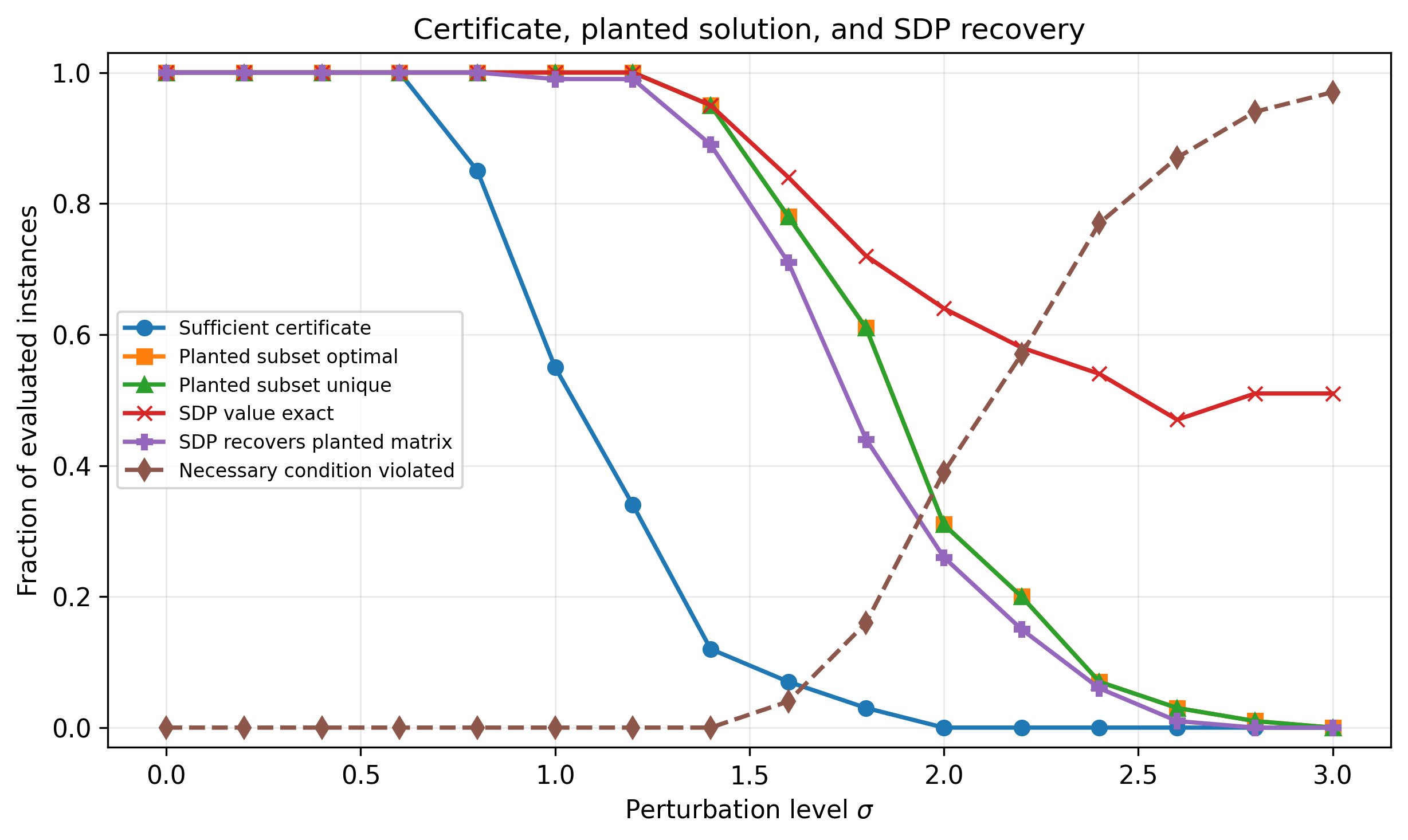}
    \caption{Empirical frequencies along the additive-perturbation paths for
\(n=50\), \(k=5\), and \(d=10\).  The curves report satisfaction of the
sufficient condition, optimality and uniqueness of the planted subset,
SDP value exactness, recovery of the planted rank-one matrix, and violation
of the necessary condition.}
    \label{fig:tightness_probabilities}
\end{figure}

\begin{table}[!ht]
\centering
\caption{Empirical frequencies at selected noise levels for \(n=50\),
\(k=5\), and \(d=10\). Under the tightness conditions, ``Guaranteed''
means that the sufficient condition holds, whereas ``Excluded'' means that
the necessary condition is violated for the planted solution \(S\).}
\label{tab:hollow-n50-k5-rates}
\footnotesize
\setlength{\tabcolsep}{5pt}
\renewcommand{\arraystretch}{1.08}
\begin{tabular}{@{}c cc cc c@{}}
\toprule
& \multicolumn{2}{c}{Tightness conditions}
& \multicolumn{2}{c}{SDP} & \\
\cmidrule(lr){2-3}
\cmidrule(lr){4-5}
\(\noise\)
& Guaranteed
& Excluded
& Exact value
& Recovers \(S\)
& \(S\) optimal \\
\midrule
0.0 & 100\% & 0\%  & 100\% & 100\% & 100\% \\
0.8 & 85\%  & 0\%  & 100\% & 100\% & 100\% \\
1.2 & 34\%  & 0\%  & 100\% & 99\%  & 100\% \\
1.8 & 3\%   & 16\% & 72\%  & 44\%  & 61\%  \\
2.2 & 0\%   & 57\% & 58\%  & 15\%  & 20\%  \\
2.6 & 0\%   & 87\% & 47\%  & 1\%   & 3\%   \\
3.0 & 0\%   & 97\% & 51\%  & 0\%   & 0\%   \\
\bottomrule
\end{tabular}
\end{table}

Figure~\ref{fig:tightness_probabilities} and
Table~\ref{tab:hollow-n50-k5-rates} summarize the empirical robustness of the
tightness conditions, the SDP relaxation, and the planted solution under
increasing perturbations. At zero noise, every reference instance satisfies the sufficient condition and has \(S\) as its unique global optimum.

The sufficient condition holds for all instances through \(\noise=0.6\),
for \(85\%\) at \(\noise=0.8\), and for \(55\%\) at \(\noise=1\), where the
root-mean-square perturbation equals the initial planted-cloud radius.
No instance remains certified at \(\noise=2\).
The SDP relaxation remains exact well beyond the certified range:
at \(\noise=1.6\), only \(7\%\) of the instances remain certified, whereas
\(84\%\) are value-exact, \(78\%\) retain the planted optimum, and \(71\%\)
recover its rank-one matrix. Certification is lost before planted optimality
on all 100 perturbation paths and before SDP value exactness on all 89 paths
on which a loss of value exactness is observed. The median first-loss noise
levels are \(1.2\) for certification and \(2.0\) for both SDP value exactness
and planted optimality. These results show that the sufficient condition
tolerates substantial perturbations, although it can be conservative
relative to the empirical recovery range.

The planted subset first loses optimality at \(\noise=1.4\) and remains
optimal in \(31\%\) of the instances at \(\noise=2\). The planted-optimality
and planted-uniqueness curves coincide throughout this run because every
observed globally optimal planted subset has a strictly worse best competitor.
This absence of ties is expected for continuously distributed random data.
SDP value exactness persists even after the planted subset ceases to be
optimal. At \(\noise=2.4\), the SDP is value-exact in \(54\%\) of the
instances, whereas \(S\) is optimal in \(7\%\). At \(\noise=3\), these
frequencies are \(51\%\) and \(0\%\), respectively, and no instance recovers
the planted matrix. Hence, late SDP value exactness reflects tightness of the
optimal value after the planted subset is no longer optimal, rather than
recovery of the prescribed planted matrix.

Finally, the necessary condition is less sensitive to increasing noise than
the sufficient condition. Its violation rate is \(87\%\) at \(\noise=2.6\),
\(94\%\) at \(\noise=2.8\), and \(97\%\) at \(\noise=3.0\). Across the 100
paths, the first necessary-condition violation occurs after the first loss of
planted optimality on 80 paths and at the same sampled noise level on 17
paths. On the remaining three paths, planted optimality is lost but the
necessary condition is not violated by the end of the grid. Thus, no
necessary-condition violation precedes the loss of planted optimality.
This delayed response is consistent with the condition ruling out SDP
tightness at the prescribed matrix rather than directly characterizing
discrete optimality. In these experiments, its violation nevertheless
provides a reliable but conservative indication of lost planted optimality.

We repeated the experiment for $k\in\{10,20\}$ and observed the same
qualitative behavior of the tightness conditions and SDP exactness.
These additional results are omitted for brevity.

\subsection{Branch-and-cut results}
\label{sub:bb_res}
\paragraph{Implementation details}
We obtain incumbent solutions with the GRASP heuristic proposed in \cite{ettahri2027grasp}. It starts from a seed point
and repeatedly adds a point sampled from a restricted candidate list whose
members have the smallest marginal sums of squared distances to the current
subset. Once the subset contains $k$ points, a best-improvement one-swap
local search exchanges one selected and one unselected point until no
improving swap remains.
Before starting the branch-and-cut algorithm, we perform 50 GRASP random starts and use the best solution to initialize the global upper bound.

As for the cutting-plane setting, we randomly separate at most $100,000$ valid cuts, sort them in decreasing order with respect to the violation, and add the first $5,000$ violated ones in each iteration. The tolerance for checking the violation is set to $10^{-4}$. 
We stop the cutting-plane procedure when the separation routine does not find violated inequalities or when the relative difference of the lower bound between two consecutive iterations is less than or equal to $10^{-3}$. We branch on the most fractional diagonal entry, i.e., the $Z_{ii}$ closest to $0.5$. The inequalities generated at each node are inherited by its children. Finally, we visit the tree with the best-first search strategy and require an optimality tolerance of $\varepsilon=10^{-4}$ on the gap, i.e., we terminate the branch-and-cut method when $(UB-LB)/\max\{1,|UB|\}\leq\varepsilon$, where UB and LB denote the incumbent objective and the best bound, respectively.

Our algorithm, named \texttt{MVkSP-SDP}, is implemented in C++ with some routines written in MATLAB (version 2022b). The SDP relaxation at each node is solved by means of SDPNAL+, a MATLAB software that implements an augmented Lagrangian method for SDPs with bound constraints \cite{sun2020sdpnal+}. We set the accuracy tolerance of SDPNAL+ to $10^{-5}$ in the relative KKT residual. We use Gurobi (version 13.0) to solve the LPs within the post-processing procedure. We run the experiments on a laptop with Intel(R) i7-13700H CPU, 32 GB of RAM, and Ubuntu 22.04.3 LTS. The source code is available at \url{https://github.com/antoniosudoso/MVkSP-SDP}.

\begin{table}[!ht]
\centering
\scriptsize
\setlength{\tabcolsep}{2.5pt}
\renewcommand{\arraystretch}{0.88}
\caption{Overview of the datasets used in the experiments.
For each dataset, ID is a sequential identifier, $n$ is the number of
{observations}, and $d$ is the number of features.}
\label{tab:datasets}
\begin{minipage}{0.48\textwidth}
\centering
\begin{tabular}{c l r r}
\hline
ID & Dataset & $n$ & $d$ \\
\hline
1  & ruspini & 75  & 2 \\
2  & caesarian & 80  & 5 \\
3  & ceramic & 88  & 17 \\
4  & immunotherapy & 90  & 7 \\
5  & analcatdata\_chlamydia & 100 & 17 \\
6  & analcatdata\_creditscore & 100 & 7 \\
7  & analcatdata\_neavote & 100 & 4 \\
8  & blogger & 100 & 15 \\
9  & fertility & 100 & 9 \\
10 & zoo & 101 & 16 \\
11 & breast-tissue & 106 & 9 \\
12 & molecular-biology & 106 & 228 \\
13 & analcatdata\_michiganacc & 108 & 17 \\
14 & cloud & 108 & 8 \\
15 & visualizing\_environmental & 111 & 3 \\
16 & robot-failures & 117 & 90 \\
17 & acute-inflammations & 120 & 11 \\
18 & rabe\_266 & 120 & 2 \\
19 & KungChi & 123 & 39 \\
20 & heart-switzerland & 123 & 12 \\
21 & fruitfly & 125 & 4 \\
22 & datatrieve & 130 & 8 \\
23 & visualizing\_livestock & 130 & 31 \\
24 & analcatdata\_seropositive & 132 & 5 \\
25 & hayes-roth & 132 & 4 \\
26 & veteran & 137 & 7 \\
27 & newton\_hema & 140 & 2 \\
28 & sleuth\_case2002 & 147 & 6 \\
29 & lymph & 148 & 47 \\
30 & iris & 150 & 4 \\
31 & teachingAssistant & 151 & 5 \\
32 & grub-damage & 155 & 32 \\
33 & auto\_price & 159 & 15 \\
34 & analcatdata\_wildcat & 163 & 5 \\
35 & wine & 178 & 13 \\
36 & backache & 180 & 32 \\
37 & triazines & 186 & 58 \\
38 & lowbwt & 189 & 9 \\
39 & wisconsin & 194 & 32 \\
40 & parkinsons & 195 & 22 \\
41 & heart-long-beach & 200 & 13 \\
42 & pwLinear & 200 & 10 \\
43 & sonar & 208 & 60 \\
44 & seismic-bumps & 210 & 7 \\
45 & new-thyroid & 215 & 5 \\
46 & chscase\_geyser & 222 & 2 \\
47 & chatfield & 235 & 12 \\
48 & qualitative-bankruptcy & 250 & 6 \\
49 & bodyfat & 252 & 14 \\
\hline
\end{tabular}
\end{minipage}
\hfill
\begin{minipage}{0.48\textwidth}
\centering
\begin{tabular}{c l r r}
\hline
ID & Dataset & $n$ & $d$ \\
\hline
50 & spect & 267 & 22 \\
51 & heart-statlog & 270 & 13 \\
52 & mu284 & 284 & 9 \\
53 & badges & 294 & 10 \\
54 & heart-h & 294 & 13 \\
55 & cleve & 303 & 13 \\
56 & cleveland-nominal & 303 & 8 \\
57 & haberman & 306 & 3 \\
58 & vertebra-column & 310 & 6 \\
59 & plasma\_retinol & 315 & 18 \\
60 & solar-flare & 315 & 12 \\
61 & pc1\_req & 320 & 8 \\
62 & visualizing\_galaxy & 323 & 4 \\
63 & ecoli & 336 & 7 \\
64 & liver-disorders & 345 & 7 \\
65 & Fish & 350 & 463 \\
66 & ionosphere & 351 & 33 \\
67 & jEdit\_4.2\_4.3 & 369 & 8 \\
68 & vinnie & 380 & 2 \\
69 & cars & 392 & 7 \\
70 & calendarDOW & 399 & 32 \\
71 & MusicEmotion & 400 & 50 \\
72 & chscase\_census & 400 & 7 \\
73 & wholesale-customers & 440 & 7 \\
74 & GunPointAgeSpan & 451 & 150 \\
75 & sa-heart & 462 & 9 \\
76 & thoracic-surgery & 470 & 16 \\
77 & clean & 476 & 166 \\
78 & PopularKids & 478 & 32 \\
79 & Computers & 500 & 720 \\
80 & ECG5000 & 500 & 140 \\
81 & EthanolLevel & 500 & 1751 \\
82 & synthetic\_control & 600 & 60 \\
83 & AbnormalHeartbeat & 606 & 3053 \\
84 & Strawberry & 613 & 235 \\
85 & ScreenType & 750 & 720 \\
86 & energy & 768 & 16 \\
87 & gene & 801 & 20531 \\
88 & FordB & 810 & 500 \\
89 & vehicle\_silhouettes & 846 & 18 \\
90 & ECGFiveDays & 884 & 136 \\
91 & UWaveGestureLibrary & 896 & 945 \\
92 & qsar\_fish & 908 & 7 \\
93 & CBF & 930 & 128 \\
94 & TwoPatterns & 1000 & 128 \\
95 & Wafer & 1000 & 152 \\
96 & banknote\_authentication & 1372 & 4\\
97 & yeast & 1484 & 8\\
98 & aspirin & 1500 & 63\\
\hline
\end{tabular}
\end{minipage}
\end{table}

\medskip
We evaluate the computational performance of the proposed \texttt{MVkSP-SDP} solver on a benchmark set composed of 98 {public datasets} available in UCI\footnote{\url{https://archive.ics.uci.edu/}} and
UCR\footnote{\url{https://www.cs.ucr.edu/~eamonn/time_series_data_2018/}}
repositories. The instances cover a broad range of problem sizes, with the number of data points varying within $n \in [75,1500]$, and the dimension ranging within $d \in [2,20531]$. An overview of the considered datasets is reported in Table~\ref{tab:datasets}.

To assess both robustness and scalability, the experiments were conducted for different values of the cardinality parameter,
$k \in \{5,10,20,30\}$, resulting in a total of 392 tested instances.
The proposed \texttt{MVkSP-SDP} solver is compared against two state-of-the-art general-purpose optimization solvers: the commercial solver \texttt{GUROBI} (version 13.0) and the open-source solver \texttt{SCIP} (version 10.0). Both solvers were applied to the nonconvex formulation~\eqref{eq:quadratic_distance}. We also tested \texttt{BARON} (version 25.3.19); however, preliminary experiments indicated that its performance was consistently inferior to that of \texttt{GUROBI} and \texttt{SCIP}, both in terms of runtime and solution quality. For this reason, \texttt{BARON} is omitted from the final comparison. Moreover, we do not include \texttt{BiqCrunch}~\cite{krislock2016computational}
or \texttt{BiqBin}~\cite{gusmeroli2022biqbin} in the comparison,
although both have been applied to H$k$SP or D$k$SP.
Their documented input formats require integer quadratic coefficients or integer edge weights, whereas generic MV$k$SP instances involve real-valued squared Euclidean distances. Rational distances can be scaled to integers
exactly, but the resulting coefficients may be prohibitively large;
rounding instead changes the problem and may alter its optimal solutions.
All experiments were executed with a time limit of 10,800 seconds (3 hours). For all solvers, the relative optimality gap tolerance was set to $10^{-4}$, consistently with the stopping criterion adopted for \texttt{MVkSP-SDP}.

\begin{figure}[!ht]
    \centering
    \includegraphics[scale=0.6]{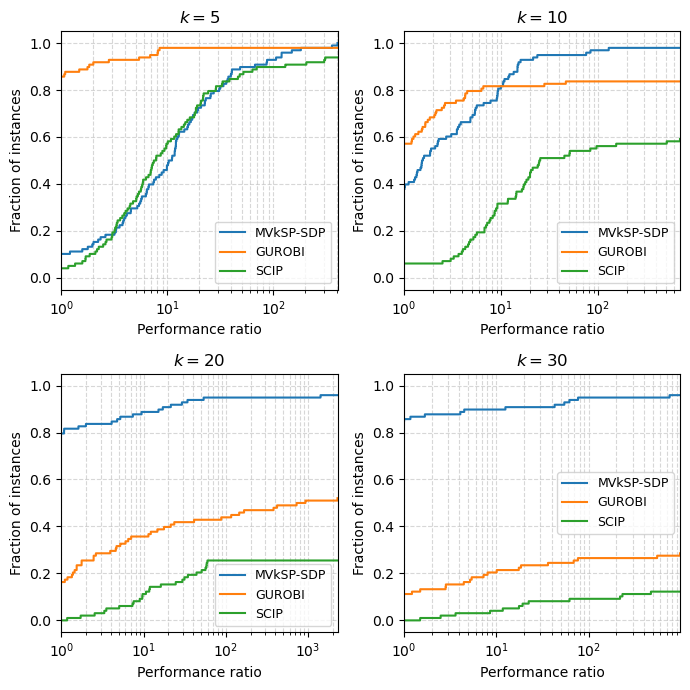}%
\caption{Performance profiles comparing \texttt{MVkSP-SDP}, \texttt{GUROBI}, and \texttt{SCIP} for $k \in \{5,10,20,30\}$.
The horizontal axis reports the performance ratio on a logarithmic scale, while the vertical axis represents the fraction of instances solved within that ratio of the best running time for each instance.}
    \label{fig:pp}
\end{figure}

Figure~\ref{fig:pp} reports the performance profiles \cite{dolan2002benchmarking} for the considered solvers and for each value of $k$. At a given point on the horizontal axis, each profile reports the fraction
of instances solved within that factor of the best runtime attained by any
solver. Higher curves indicate better performance, and their limiting heights
give the fractions of instances solved within the time limit.
The profiles show a clear change in relative performance as $k$ increases.
\texttt{GUROBI} has the best overall profile for $k=5$, while \texttt{MV$k$SP-SDP}
generally outperforms \texttt{SCIP}. For $k=10$, \texttt{MVkSP-SDP} leads over most
performance ratios, and its advantage becomes more pronounced for
$k=20$ and $k=30$. In contrast, the profiles of \texttt{GUROBI} and especially \texttt{SCIP} flatten considerably as $k$ increases, reflecting the growing number of instances reaching the time limit. 

Figure~\ref{fig:ppn} in Appendix \ref{appendix:size_profiles} further analyzes the performance profiles by separating the instances according to their size, distinguishing between small instances ($n < 300$) and large instances ($n \geq 300$). For $k=5$, \texttt{GUROBI} remains competitive on both small and large instances, although \texttt{MVkSP-SDP} becomes increasingly competitive as the instance size grows. 
For $k=10$, \texttt{MVkSP-SDP} already dominates on large instances, while \texttt{GUROBI} retains an advantage on smaller ones. For $k=20$ and $k=30$, \texttt{MVkSP-SDP} clearly outperforms both competing solvers on large-scale instances, solving the vast majority of the problems within small performance ratios, whereas \texttt{GUROBI} and \texttt{SCIP} experience a degradation in performance and frequently reach the time limit. 
Together, these results indicate that the proposed
method scales more effectively with both $n$ and $k$.

Table~\ref{tab:comparison_by_k} complements the runtime profiles with completion rates and optimality gaps. Across all 392 instances, \texttt{MVkSP-SDP} reaches the
time limit on 10 instances, compared with 135 for \texttt{GUROBI} and 205 for \texttt{SCIP}. Although \texttt{GUROBI} has the lowest average runtime for $k=5$, \texttt{MVkSP-SDP} solves every instance in this group and has the lowest average runtime for all larger values of $k$. At $k=30$, its timeout rate is 4.08\%, compared with 71.43\% for \texttt{GUROBI} and 87.76\% for \texttt{SCIP}. The residual gaps further distinguish the methods. On the instances reaching the time limit, the average gap of \texttt{MVkSP-SDP} ranges from 0.11\% to 0.39\%, and its maximum gap remains below 1\%.
In comparison, the average gaps of \texttt{GUROBI} and \texttt{SCIP} exceed 64\% and 93\%, respectively, for $k\geq10$.
Thus, even when \texttt{MVkSP-SDP} does not terminate within the imposed time limit, it is still able to produce solutions very close to optimality.

\begin{table}[!ht]
\centering
\scriptsize
\setlength{\tabcolsep}{3pt}
\caption{Performance comparison by value of $k$. The table reports the percentage of instances reaching the time limit (TL), the average runtime in seconds, the average optimality gap, and the average and maximum optimality gaps among the instances reaching the time limit.}
\label{tab:comparison_by_k}
\begin{tabular}{rlrrrrr}
\toprule
$k$ & Solver & TL (\%) & \shortstack{Avg. time\\(s)} & \shortstack{Avg. gap\\(\%)} & \shortstack{Avg. TL gap\\(\%)} & \shortstack{Max. TL gap\\(\%)} \\
\midrule
5  & \texttt{MVkSP-SDP} & 0.00 & 446.87  & 0.00  & --       & --     \\
5  & \texttt{GUROBI}    & 2.04 & 287.10  & 1.49  & 72.44    & 96.12  \\
5  & \texttt{SCIP}      & 6.12 & 898.24  & 6.12  & 100.00   & 100.00 \\
\midrule
10 & \texttt{MVkSP-SDP} & 2.04 & 572.89  & 0.00  & 0.11    & 0.16   \\
10 & \texttt{GUROBI}    & 16.33 & 1885.63 & 10.50 & 64.31  & 99.07  \\
10 & \texttt{SCIP}      & 40.82 & 4649.08 & 38.15 & 93.47  & 100.00 \\
\midrule
20 & \texttt{MVkSP-SDP} & 4.08 & 875.80  & 0.02  & 0.39   & 0.90   \\
20 & \texttt{GUROBI}    & 47.96 & 5736.25 & 31.39 & 65.46  & 99.93  \\
20 & \texttt{SCIP}      & 74.49 & 8136.96 & 72.11 & 96.81  & 100.00 \\
\midrule
30 & \texttt{MVkSP-SDP} & 4.08 & 1266.87 & 0.01  & 0.29   & 0.94   \\
30 & \texttt{GUROBI}    & 71.43 & 7928.24 & 47.69 & 66.77  & 100.00 \\
30 & \texttt{SCIP}      & 87.76 & 9572.59 & 85.46 & 97.39  & 100.00 \\
\bottomrule
\end{tabular}
\end{table}

Finally, we connect the computational performance of \texttt{MVkSP-SDP} with the
tightness analysis of Section~\ref{sec:tightness} by examining the instances solved by the
initial root relaxation. The branch-and-cut algorithm proves optimality before
adding cutting planes or branching on 61 of the 392 benchmark instances
(15.56\%). Applying the certificate of Algorithm~\ref{alg:exact-recovery} to the optimal subsets
obtained in these instances certifies exact recovery in 50 cases (81.97\%). For the
remaining instances, global optimality cannot be certified through this
condition and must instead rely on the exact solution framework, using
cutting planes and branching whenever required.

\section{Conclusions}\label{sec:concl}
We studied the minimum variance \(k\)-subset problem using semidefinite
programming.  By identifying the face containing the feasible region, we
obtained an equivalent facially reduced relaxation satisfying Slater's
condition and strong duality. The resulting primal--dual characterization led to directly checkable geometric conditions for SDP tightness and exact recovery.  These conditions capture the interaction between the internal shape of a candidate subset, its separation from the unselected points, and the separation among the unselected points.
The experiments support the practical relevance of these
conditions.  The sufficient condition guarantees exact recovery over a substantial perturbation range, although the SDP in practice frequently remains exact beyond this range.  At higher noise levels, violation of the necessary condition provides a conservative indication of lost planted optimality.
We also embedded the SDP relaxation into a branch-and-cut
method.  The resulting algorithm is particularly effective as \(n\) and
\(k\) increase, substantially outperforming the general-purpose solvers on
the most challenging instances with up to \(1{,}500\) observations.

Several directions remain open.  From a theoretical perspective, strengthening
the sufficient condition could narrow the gap between the guaranteed and
empirically exact regimes.  A probabilistic analysis of recovery thresholds
under suitable random models could further explain the transitions observed
in the numerical experiments. From an algorithmic perspective, a tailored solver that exploits the low-rank structure of the SDP relaxation could improve scalability to larger instances.  

\section*{Statements and declarations}
\paragraph{Competing interests.}
The authors declare that they have no competing interests.

\paragraph{Data and code availability.}
The benchmark datasets are publicly available. The source code is available at \url{https://github.com/antoniosudoso/MVkSP-SDP}.

%
%
%





\bibliographystyle{spmpsci}
\bibliography{abbr, biblio} 




\begin{appendix}

\section{Proof of Lemma \ref{lemma:centering_identities}}
\label{appendix:centering_identities}
\begin{proof}
Using the identity $\sum_{j \in S} r_j = 0$, we have
\begin{align*}
\sum_{j\in S}D_{ij}
&=
\sum_{j\in S}\|r_i-r_j\|^2
=
k\|r_i\|^2+\sum_{j\in S}\|r_j\|^2 - 2 r_i^\top \left(\sum_{j \in S} r_j \right) \\
&= kq_i+\sum_{j\in S}q_j
- 2 r_i^\top \left(\sum_{j \in S} r_j \right)
=
k(q_i + \bar{q}),
\end{align*}
which proves \eqref{eq:centering_Dsum}. Next, from $v^\top \one_k = 0$ we have $\sum_{i,j\in S}v_iv_j
= 0$.
Thus, using $D_{ij} = q_i + q_j - 2 r_i^\top r_j$, the two radial terms vanish
after summation and
\[
v^\top D_{SS}v
=
-2\sum_{i,j\in S}v_iv_jr_i^\top r_j = -2 \left\|
\sum_{i\in S}v_ir_i
\right\|^2
=
-2 \|R_Sv\|^2,
\]
which proves \eqref{eq:centering_DSS}.
\qed
\end{proof}

\section{Proof of Theorem \ref{thm:spherical-tightness}}
\label{appendix:proof_spherical}
\begin{proof}
We construct a dual certificate $(u,\tau)$ satisfying the three conditions of Proposition~\ref{prop:tightness}. Define
\begin{equation*}
  \gamma \;:=\; \min\!\left\{
      \frac{\Delta_\times-\Delta_S}{2k},\;
      \frac{\Delta_{\bar S}-\Delta_S}{4k}
  \right\}, \quad \alpha \;:=\; -\frac{\Delta_S}{2k},
\end{equation*}
and set 
\begin{equation*}
    u_i = \alpha \quad \forall i \in S, \qquad u_j = \alpha - \gamma \quad \forall j \in \bar{S}, \qquad \tau = \frac{\bar q}{k}.
\end{equation*}
Recall that,
\[
H_{ij}
=
\frac{D_{ij}}{2k^2}
-u_i\delta_{ij}
+\frac{u_i+u_j}{2k}
-\tau\delta_{ij},
\]
where $\delta_{ij}$ is the Kronecker delta. We first prove condition $H_{SS}\mathbf{1}_k = \mathbf{0}_k$. 
Since $u_i=\alpha$ for all $i\in S$, the terms involving $\alpha$ cancel in
the row sum, and therefore
\[
(H_{SS}\mathbf 1_k)_i
=
\frac1{2k^2}\sum_{j\in S}D_{ij}-\tau.
\]
Now, using \eqref{eq:centering_Dsum} from Lemma~\ref{lemma:centering_identities} and the spherical assumption $q_i = \bar q$ for all $i \in S$, we have
\[
(H_{SS}\mathbf 1_k)_i
=
\frac1{2k}(q_i + \bar q) - \frac{\bar q}{k} = 0,
\]
which proves condition \eqref{eq:block_kernel}.

We now prove that $H_{SS}\succeq0$. Since $H_{SS}\mathbf 1_k=\mathbf 0_k$,
it suffices to show that $v^\top H_{SS}v\ge0$ for all $v\in\mathbb R^k$ such that 
$v^\top\mathbf 1_k=0$. Using the definition of $H_{SS}$ together with $u_i=\alpha$ for all $i \in S$,
we obtain
\[
v^\top H_{SS}v
=
\frac1{2k^2}v^\top D_{SS}v
-\alpha\|v\|^2
-\tau\|v\|^2,
\]
because the term $\frac{\alpha}{k}\mathbf 1_k\mathbf 1_k^\top$
vanishes on vectors orthogonal to $\mathbf 1_k$.
Substituting~\eqref{eq:centering_DSS} from Lemma~\ref{lemma:centering_identities} and using $\alpha = -\Delta_S/(2k)$ and {$\tau = \bar q/k$}
gives
\begin{equation}
v^\top H_{SS}v
=
-\frac1{k^2}\|R_Sv\|^2
+
\left(
\frac{\Delta_S}{2k}-\frac{\bar q}{k}
\right)\|v\|^2.
\label{eq:H-quad}
\end{equation}
It therefore remains to bound $\|R_Sv\|^2$.
Let $G:=R_S^\top R_S$ be the Gram matrix of the centered points. Since $D_{ij}=2\bar{q}-2G_{ij}$
and $D_{ij}\le\Delta_S$ for all $i,j\in S$, we have $G_{ij}\ge \bar{q}-\frac{\Delta_S}{2}$. Define {$c:={\Delta_S}/2-\bar q\ge0$},
and set $A:=G+c\,\mathbf 1_k\mathbf 1_k^\top$.
The nonnegativity of $c$ follows from
$2k\bar q=\sum_{j\in S}D_{ij}\le (k-1)\Delta_S$ for any $i\in S$.
Then, $A$ is entrywise nonnegative. Moreover, $A\mathbf 1_k
=
G\one_k
+
ck\one_k
=
ck \one_k$,
since $G\mathbf 1_k=\zero_k$.
Thus every row of $A$ sums to $ck$. Since $A$ is symmetric and entrywise
nonnegative, $x^\top Ax
\le
ck\|x\|^2$ for every $x\in\mathbb R^k$.
Now, because $v^\top\one_k=0$, the equality $v^\top Gv
=
v^\top Av$ holds.
Therefore,
\[
\|R_Sv\|^2
=
v^\top Gv
=
v^\top Av
\le
ck\,\|v\|^2
=
k\left(
\frac{\Delta_S}{2}-\bar q
\right)\|v\|^2.
\]
Substituting this estimate into~\eqref{eq:H-quad} yields
\[
v^\top H_{SS}v
\ge
-\frac1{k^2}
k\left(
\frac{\Delta_S}{2}-\bar{q}
\right)\|v\|^2
+
\left(
\frac{\Delta_S}{2k}-\frac{\bar{q}}{k}
\right)\|v\|^2
=
0.
\]
Hence $v^\top H_{SS}v\ge0$ for all $v\perp\mathbf 1_k$.
Consequently, condition \eqref{eq:block_psd} holds.

Finally, we verify $H_{ij}>0$ for all $(i,j)\notin S\times S$ by distinguishing three cases.

\noindent
\emph{Case~1: $i\in S$, $j\in\bar S$.}
With $u_i=\alpha$ and $u_j=\alpha-\gamma$,
\[
  H_{ij}
  = \frac{D_{ij}}{2k^2}+\frac{u_i+u_j}{2k}
  = \frac{D_{ij}-\Delta_S}{2k^2} - \frac{\gamma}{2k}.
\]
Since $D_{ij}\ge \Delta_\times$ and $\gamma\le(\Delta_\times-\Delta_S)/(2k)$,
\[
  H_{ij}
  \ge \frac{\Delta_\times-\Delta_S}{2k^2} - \frac{\Delta_\times-\Delta_S}{4k^2}
  = \frac{\Delta_\times-\Delta_S}{4k^2} > 0.
\]
\noindent
\emph{Case~2: $i,j\in\bar S$, $i\neq j$.}
With $u_i=u_j=\alpha-\gamma$,
\[
  H_{ij}
  = \frac{D_{ij}}{2k^2}+\frac{\alpha-\gamma}{k}
  = \frac{D_{ij}-\Delta_S}{2k^2} - \frac{\gamma}{k}.
\]
Since $D_{ij}\ge \Delta_{\bar S}$ and $\gamma\le(\Delta_{\bar S}-\Delta_S)/(4k)$,
\[
  H_{ij}
  \ge \frac{\Delta_{\bar S}-\Delta_S}{2k^2} - \frac{\Delta_{\bar S}-\Delta_S}{4k^2}
  = \frac{\Delta_{\bar S}-\Delta_S}{4k^2} > 0.
\]

\noindent
\emph{Case~3: $j\in\bar S$.}
With $u_j=\alpha-\gamma$,
\[
  H_{jj}
  = -u_j\frac{k-1}{k} - \tau
  = \Bigl(\frac{\Delta_S}{2k}+\gamma\Bigr)\frac{k-1}{k} - \frac{\bar q}{k}.
\]

\noindent
Note that {$D_{ij}\le \Delta_S$} for all $i,j\in S$. Fix any $i\in S$. Then
\[
\sum_{j\in S}D_{ij}
=
D_{ii}+\sum_{\substack{j\in S\\ j\neq i}}D_{ij}
=
\sum_{\substack{j\in S\\ j\neq i}}D_{ij} \le (k-1)\Delta_S,
\]
since the last sum contains $k-1$ terms and each of them is at most $\Delta_S$.
Using the identity \eqref{eq:centering_Dsum} in Lemma \ref{lemma:centering_identities}, we obtain
$2k\bar q\le (k-1)\Delta_S$,
or equivalently, $\Delta_S\ge \frac{2k\bar q}{k-1}$.
Therefore,
\[
  H_{jj}
  \ge \left(\frac{2k\bar q}{2k(k-1)} + \gamma \right) \frac{k-1}{k} - \frac{\bar q}{k}
  = \gamma\frac{k-1}{k} > 0.
\]
Thus condition \eqref{eq:block_positive} holds. All three conditions of Proposition~\ref{prop:tightness} are satisfied, so
$Z^S=x^S(x^S)^\top$ is the unique optimal solution of
the SDP relaxation.
\qed
\end{proof}

\section{Proof of Lemma~\ref{lemma:propagation}}
\label{appendix:lemma_propagation}
\begin{proof}
Since $Z\succeq 0$, there exist vectors $w_1,\dots,w_n$ such that $Z_{ij}=w_i^\top w_j$ for all $i, j \in \{1, \dots, n\}$. Let $\bar{w}=\sum_{j=1}^n w_j$. Using the constraint $Z\mathbf{1}_n=k\,\operatorname{diag}(Z)$, we obtain for each $i \in \{1, \dots, n\}$,
\[
w_i^\top \bar{w} = \sum_{j=1}^n w_i^\top w_j = (Z\mathbf{1}_n)_i = k Z_{ii} = k\|w_i\|^2.
\]
Moreover,
\[
\|\bar w\|^2 = \sum_{i=1}^n \sum_{j=1}^n w_i^\top w_j = \mathbf{1}_n^\top Z \mathbf{1}_n.
\]
Since $Z\mathbf{1}_n = k\,\operatorname{diag}(Z)$ and $\operatorname{tr}(Z)=k$, it follows that $\mathbf{1}_n^\top Z \mathbf{1}_n = k \operatorname{tr}(Z)=k^2$, hence $\|\bar{w}\|=k$. Applying the Cauchy--Schwarz inequality gives $k\|w_i\|^2
= w_i^\top \bar w
\leq \|w_i\|\,\|\bar w\|
= k\|w_i\|$, which implies $\|w_i\|\le 1$, that is, $0\le Z_{ii}=\|w_i\|^2 \le 1$.
If $Z_{ii}=0$, then $\|w_i\|=0$, so $w_i=0$. Therefore $Z_{ij}=w_i^\top w_j=0$ for all $j \in \{1, \dots, n\}$.

If instead $Z_{ii}=1$, then $\|w_i\|=1$, and the inequality
$w_i^\top \bar w \le \|w_i\| \|\bar w\|$ holds with equality, since $w_i^\top \bar{w} = k\|w_i\|^2 = k = \|w_i\| \|\bar{w}\|$.
Hence equality holds in Cauchy--Schwarz, and therefore $w_i$ and $\bar{w}$ are collinear. Since $\|w_i\|=1$ and $\|\bar{w}\|=k$, we must have $w_i=\frac{\bar{w}}{k}$.
Consequently, for every $j \in \{1, \dots, n\}$,
\[
Z_{ij}=w_i^\top w_j=\frac{1}{k}\bar{w}^\top w_j
=\frac{1}{k}\sum_{\ell=1}^n w_\ell^\top w_j
=\frac{1}{k}(Z\mathbf{1}_n)_j
=Z_{jj}.
\]
Finally, consider the principal submatrix $\widehat Z$ obtained by deleting the $i$-th row and column of $Z$. In the case $Z_{ii}=0$, the deleted row and column are zero, so $\widehat Z$ remains feasible for the SDP on the remaining $n-1$ indices with the same cardinality parameter $k$. In the case $Z_{ii}=1$, we have, for every $j\neq i$,
\[
\sum_{\ell\neq i} \widehat Z_{j\ell}
=
\sum_{\ell=1}^n Z_{j\ell}-Z_{ji}
=
kZ_{jj}-Z_{jj}
=
(k-1)Z_{jj},
\]
and also
\[
\operatorname{tr}(\widehat Z)=\operatorname{tr}(Z)-Z_{ii}=k-1.
\]
Thus $\widehat Z$ satisfies the same type of constraints with cardinality parameter $k-1$. Since principal submatrices of positive semidefinite matrices preserve these properties, the reduced matrix is feasible for the corresponding SDP on $n-1$ indices.
\qed
\end{proof}

\clearpage

\section{Performance profiles by instance size}
\label{appendix:size_profiles}
\begin{figure}[!ht]
    \centering
    \includegraphics[scale=0.54]{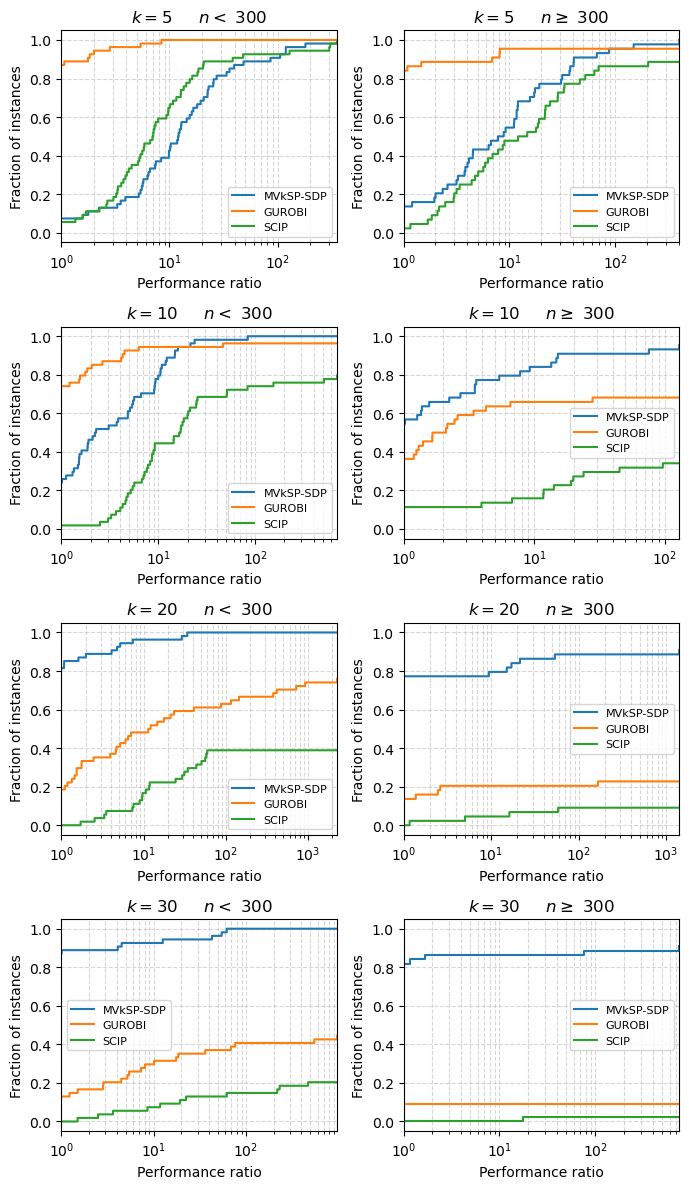}%
\caption{Performance profiles comparing \texttt{MVkSP-SDP}, \texttt{GUROBI}, and \texttt{SCIP} for $k \in \{5,10,20,30\}$.
For each value of $k$, instances are partitioned according to the problem size $n$ into small instances ($n < 300$) and large instances ($n \geq 300$). The horizontal axis reports the performance ratio on a logarithmic scale.}
    \label{fig:ppn}
\end{figure}

\end{appendix}

\end{document}